\documentclass[a4paper]{article}

\usepackage[leqno,intlimits]{amsmath}
\usepackage{amssymb}
\usepackage{amsthm}
\usepackage{hyperref}
\usepackage[usenames,dvipsnames]{pstricks}
\usepackage{epsfig}
\usepackage{appendix}

\newtheorem{tm}{Theorem}[section]

\newtheorem{rem}[tm]{Remark}

\newcommand*{\un}[1]{\underline{#1}}
\newcommand*{\Zb}{\mathbb Z}
\newcommand*{\Rb}{\mathbb R}
\newcommand*{\om}{\omega}
\newcommand*{\de}{\delta}
\newcommand*{\ze}{\zeta}
\newcommand*{\vp}{\varphi}
\newcommand*{\la}{\lambda}
\newcommand*{\vr}{\varrho}
\newcommand*{\Hc}{\mathcal H}
\newcommand*{\Ev}{{\bf E}}
\newcommand*{\omin}{\om^{\text{min}}}
\newcommand*{\omax}{\om^{\text{max}}}
\newcommand*{\hop}{\bigskip\noindent}
\newcommand*{\ba}{\begin{aligned}}
\newcommand*{\ea}{\end{aligned}}
\newcommand*{\be}{\begin{equation}}
\newcommand*{\ee}{\end{equation}}
\newcommand*{\te}{\theta}
\newcommand*{\si}{\sigma}
\newcommand*{\wt}{\widetilde}
\newcommand*{\up}{\uparrow}
\newcommand*{\dn}{\downarrow}
\newcommand*{\e}[1]{\text{\rm e}^{#1}}
\newcommand*{\di}{\,\text{\rm d}}

\numberwithin{equation}{section}

\begin{document}

\title{Random walk of second class particles in product shock measures}

\author{M\'arton Bal\'azs\footnote{Department of Stochastics, Institute of Mathematics, Budapest University of Technology and Economics, 1.\ Egry J\'ozsef u.\ H \'ep.\ V.\ 7., 1111 Budapest, Hungary. Email:~{\tt balazs@math.bme.hu}},
%\affil{\normalsize \it Department of Stochastics, Institute of Mathematics, Budapest University of Technology and Economics, 1 Egry J\'ozsef u.\ H V.\ 7., 1111 Budapest, Hungary}
Gy\"orgy Farkas\footnote{Budapest University of Technology and Economics, Budapest, Hungary. Email:~{\tt fgyorgy@math.bme.hu}},
P\'eter Kov\'acs\footnote{Budapest University of Technology and Economics, Budapest, Hungary. Email:~{\tt kpeter@math.bme.hu}},
%\affil{\normalsize \it Budapest University of Technology and Economics, Budapest, Hungary}
Attila R\'akos\footnote{Research Group for Condensed Matter Physics, Hungarian Academy of Sciences, Budapest University of Technology and Economics, Budafoki~u.~8, 1111 Budapest, Hungary. Email:~{\tt rakos@phy.bme.hu}}}
%\affil{\normalsize \it Research Group for Condensed Matter Physics, Hungarian Academy of Sciences, Budapest University of Technology and Economics, Budafoki~u.~8, 1111 Budapest, Hungary}
\date{\today}
\maketitle

\begin{abstract}
We consider shock measures in a class of conserving stochastic particle systems on $\Zb$. These shock measures have a product structure with a step-like density profile and include a second class particle at the shock position. We show for the asymmetric simple exclusion process, for the exponential bricklayers' process, and for a generalized zero range process, that under certain conditions these shocks, and therefore the second class particles, perform a simple random walk. Some previous results, including random walks of product shock measures and stationary shock measures seen from a second class particle, are direct consequences of our more general theorem. Multiple shocks can also be handled easily in this framework. Similar shock structure is also found in a nonconserving model, the branching coalescing random walk, where the role of the second class particle is played by the rightmost (or leftmost) particle.
\end{abstract}


\hop

\section{Introduction}

On a macroscopic level driven diffusive systems are often described by a set of conservation laws for the densities. 
These hydrodynamical equations are in general nonlinear PDEs, which can develop singularities in the solution. Shocks are discontinuities in these weak solutions, which travel with a speed also known as the Rankine-Hugoniot velocity. 
Whereas the large scale continuous description of shocks is well established \cite{Lax1973}, much less is known about the microscopic structure and dynamics, which has become a subject of intense investigation in recent years \cite{dls,ffshock,semuso,valak,qse,sokvalak,bcrw,Rakos2004,Jafarpour2004,schtaba,Jafarpour2007a,Jafarpour2008a,PreSimon2009}. 

In order to make the first steps in this direction one has to define the position of shocks on the lattice scale, which is already a nontrivial task in general. 
It is well known that second class particles, which move stochastically and follow the trajectories of density fluctuations, are attracted by shocks and therefore serve as good markers for the shock position. 
Derrida et al.\ in \cite{dls} derive the time invariant shock measure in the asymmetric simple exclusion process (ASEP) as seen from such a second class particle. 
They observe that when a certain condition holds for the asymmetry and the limiting densities, the invariant measure becomes a Bernoulli product measure with a simple step-like density profile. 
Similar product shock structure, as seen from a second class particle, was found later in another stochastic lattice model, the exponential bricklayers' process (BLP) \cite{valak}. 

A different approach, initiated by Belitsky and Sch\"utz, attempted to capture not only the structure but also the microscopic dynamics of shocks.  
In \cite{qse} they show that under the same condition as in \cite{dls} there is a family of product measures $\mu_k, k\in \Zb$ with a step-like density which evolve into linear combinations of similar measures, and the interpretation is that the shock position $k$ performs a simple random walk. 
The random walking shocks were shown to exist later in the exponential BLP \cite{sokvalak} too. 
The advantage of this description is that it doesn't use second class particles therefore it can be applied also in cases where second class particles cannot be defined (or the number of them is not conserved).
An example of such case is the branching coalescing random walk.
Although this is a non-conserving system, there are shocks with similar structure and evolution here as well \cite{bcrw}.
It is interesting that random walking shocks have also been found in systems with more than one conserved quantities \cite{Rakos2004}.

These two types of results naturally raise the question whether the second class particle itself, attracted by the shock, performs a simple random walk. 
In this paper we give an answer to this question by considering shock measures with second class particles at the shock position. 
The idea of considering such shock measures appeared in the context of the ASEP with open
boundaries in \cite{bcrw}, where a conjecture is formulated saying that the above random walk property should hold for the shock measures {\em with second class particles} too.
In a fairly general framework we show the random walk property for these measures in the ASEP, in the exponential BLP and in a generalized zero range process (GZRP) (where negative particle numbers can also occur).
While our result clearly shows the simple random walk of second class particles, the results of \cite{dls,valak,qse,sokvalak} for product shock measures also follow. Notice that the existing stationary product distribution results did not include any random walk dynamics, and the random walk results did not include the second class particle. Hence our result is genuinely new, and also connects the two types of arguments.

%This result extends previous random walk results where no information on second class particles was given.

The diffusion coefficient of a general shock in ASEP was computed and the diffusive behavior was investigated by Ferrari and Fontes \cite{ffshock}. They used the second class particle as a characterizing object for the shock location, and their result on how the shock location depends on the initial configuration of ASEP made it possible to generalize the diffusive scale-results to the case of multiple shocks with second class particles in Ferrari, Fontes and Vares \cite{semuso}.

Multiple shocks, i.e., several steps in the density profile have also been studied in \cite{qse} and \cite{sokvalak}. In this case the exact microscopic description is given, and involves several shock positions, which perform an interacting simple random walk. Due to the attraction of these ``micro-shocks'' they form a bound state with a finite width and can be considered as a single shock with a more complex structure. 
It is interesting that while the result for multiple shocks is a direct generalization of that for a single shock in the BLP \cite{sokvalak}, in the ASEP such a naive generalization does not hold. In order to be able to handle multiple shocks in the ASEP, extra particles or vacancies had to be introduced at the shock position \cite{qse}. 
Our description with the second class particles at the shock positions explains the interaction between shocks in a very natural way, without the need of artificial particles. In fact we show that previously known forms of random walking shocks can be obtained as an appropriate mixture of the two marginals of our coupling shock measure. 

It is important to note that the existence of random walking shocks is closely related to the exact solvability of these particle systems in an open geometry. 
The matrix product ansatz is a method widely used for finding stationary states of stochastic lattice models \cite{Blythe2007}. 
Recently it has been shown \cite{Jafarpour2004,Jafarpour2007a} that the occurrence of a single random walking shock implies the existence of a two-dimensional representation of the quadratic algebra appearing in the matrix product ansatz. 
Similarly, multiple shocks correspond to other finite dimensional representations \cite{Jafarpour2004,Jafarpour2007a,Jafarpour2008a}. 
Moreover, as a very recent progress, utilizing the single-particle properties of random walking shocks, the Bethe ansatz has been successfully applied \cite{PreSimon2009} in the open-boundary ASEP for the evaluation of the spectrum and also for the current large deviation function.

The paper is organized as follows. 
In section \ref{sc:fam} we introduce a wide family of stochastic particle systems on $\Zb$ with nearest neighbour jumps, and summarize their basic properties, such as stationary product measures and hydrodynamic limit. 
A few specific examples are considered in more detail. 
A surface growth interpretation can also be given to these models which is sometimes more natural. 
In section \ref{sc:results} we formulate our main results for random walking single and multiple shocks, proofs are given in section \ref{sc:proofs}. 
A different model, the branching coalescing random walk is investigated in section \ref{sc:bcrw}. This model is not in the family considered before. 
We prove a statement here which is similar in spirit to that of section \ref{sc:results}. The second class particle is replaced here by the rightmost (or leftmost) particle. 

\section{A family of models}\label{sc:fam}

The class of stochastic interacting systems we consider here appeared several times in the literature, we repeat a description recently formulated in \cite{varj2nd}. The class is a generalization of the so-called misanthrope process. We use a surface growth interpretation, but many members of this class can be understood in terms of particles jumping on the one dimensional lattice. For $-\infty\le\omin\le0$ and $1\le\omax\le\infty$ (possibly infinite valued) integers, we define the single-site state space
\[
I:\,=\left\{z\in\Zb\,:\,\omin-1<z<\omax+1\right\}
\]
and the configuration space
\[
\Omega=\left\{\un\om=(\om_i)_{i\in\Zb}\ :\ \om_i\in I\right\}=I^{\Zb}.
\]

For each pair of neighboring sites $i$ and $i+1$ of $\Zb$, we consider a column built of bricks above the edge $(i,\,i+1)$. The height of this column is denoted by $h_i$. A state configuration $\un{\om}\in\Omega$ has components $\om_i=h_{i-1}-h_i\,\in I$, being the negative discrete gradients of the height of the ``wall''. For \(\un\om\in\Omega\) and \(i\ne j\) let \(\un\om^{i,j}\) be the configuration with components
\be
\om^{i,j}_k=\left\{\ba
&\om_k,&&\text{for }k\ne i,\,j,\\
&\om_k-1,&&\text{for }k=i,\\
&\om_k+1,&&\text{for }k=j.
\ea\right.\label{eq:ugrdef}
\ee
Also, define, for a vector \(\un h\) of heights, \(\un h^{i\up}\) and \(\un h^{i\dn}\) by
\[
h^{i\up}_k=\left\{\ba
&h_k,&&\text{for }k\ne i,\\
&h_k+1,&&\text{for }k=i,
\ea\right.\qquad
h^{i\dn}_k=\left\{\ba
&h_k,&&\text{for }k\ne i,\\
&h_k-1,&&\text{for }k=i.
\ea\right.
\]

\begin{figure}
\begin{center}
\scalebox{1} % Change this value to rescale the drawing.
{
\begin{pspicture}(0,-3.32)(8.200937,3.32)
\pscircle[linewidth=0.04,dimen=outer,fillstyle=solid,fillcolor=black](3.6809375,-2.0){0.3}
\psline[linewidth=0.04](2.1809375,-2.3)(2.1809375,-2.5)(8.180938,-2.5)(8.180938,-2.3)
\psline[linewidth=0.04cm](3.1809375,-2.3)(3.1809375,-2.5)
\psline[linewidth=0.04cm](4.1809373,-2.3)(4.1809373,-2.5)
\psline[linewidth=0.04cm](5.1809373,-2.3)(5.1809373,-2.5)
\psline[linewidth=0.04cm](6.1809373,-2.3)(6.1809373,-2.5)
\psline[linewidth=0.04cm](7.1809373,-2.3)(7.1809373,-2.5)
\pscircle[linewidth=0.04,dimen=outer](5.6809373,-3.0){0.3}
\usefont{T1}{ptm}{m}{n}
\rput(3.5823438,-2.79){$i$}
\usefont{T1}{ptm}{m}{n}
\rput(4.7723436,-2.79){$i+1$}
\pscircle[linewidth=0.04,dimen=outer,fillstyle=solid,fillcolor=black](3.6809375,-1.2){0.3}
\pscircle[linewidth=0.04,dimen=outer,fillstyle=solid,fillcolor=black](6.6809373,-2.0){0.3}
\psline[linewidth=0.04cm,arrowsize=0.05291667cm 2.0,arrowlength=1.4,arrowinset=0.4]{->}(0.9809375,-3.3)(0.9809375,-0.3)
\usefont{T1}{ptm}{m}{n}
\rput(0.66234374,-1.79){$\omega$}
\psline[linewidth=0.04](2.6809375,2.9)(3.6809375,2.9)(3.6809375,0.9)(5.6809373,0.9)(5.6809373,1.9)(6.6809373,1.9)(6.6809373,0.9)(7.6809373,0.9)
\psline[linewidth=0.04,linestyle=dashed,dash=0.16cm 0.16cm](3.6809375,1.9)(4.6809373,1.9)(4.6809373,0.9)
\usefont{T1}{ptm}{m}{n}
\rput(4.182344,0.61){$i$}
\psline[linewidth=0.04cm,arrowsize=0.05291667cm 2.0,arrowlength=1.4,arrowinset=0.4]{->}(0.9809375,0.3)(0.9809375,3.3)
\usefont{T1}{ptm}{m}{n}
\rput(0.63234377,1.81){$h$}
\usefont{T1}{ptm}{m}{n}
\rput(2.7323437,-2.79){$i-1$}
\usefont{T1}{ptm}{m}{n}
\rput(3.1323438,0.61){$i-1$}
\usefont{T1}{ptm}{m}{n}
\rput(5.1723437,0.61){$i+1$}
\psarc[linewidth=0.02,arrowsize=0.05291667cm 2.0,arrowlength=1.4,arrowinset=0.4]{<-}(4.1809373,-1.1){0.6}{30.963757}{149.03624}
\psline[linewidth=0.02cm,arrowsize=0.05291667cm 2.0,arrowlength=1.4,arrowinset=0.4]{->}(4.1809373,1.0)(4.1809373,1.8)
\end{pspicture} 
}
\end{center}
\caption{Shown is the relation between the surface (above) and the particle (below) interpretation. The empty circle
denotes an antiparticle.}
\label{fig:h_omega}
\end{figure}
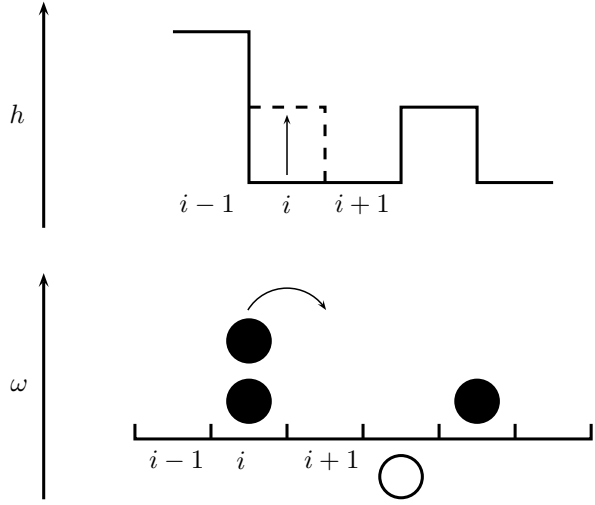

The continuous time evolution is described by jump processes. A brick can be added:
\[
\left.
\ba
\un\om&\longrightarrow\un\om^{i,i+1}\\
\un h&\longrightarrow h^{i\up}
\ea
\right\}
\text{with rate}\ p(\om_i,\,\om_{i+1}),%\label{eq:badd}
\]
or removed:
\[
\left.
\ba
\un\om&\longrightarrow\un\om^{i+1,i}\\
\un h&\longrightarrow h^{i\dn}
\ea
\right\}\text{with rate}\ q(\om_i,\,\om_{i+1}).%\label{eq:brem}
\]
Conditionally on $\un\om(t)$, these moves are independent. See Figure \ref{fig:h_omega} for visualization.
We impose the following assumptions on the rates:
\begin{itemize}
\item The rates must satisfy
\[
p(\omin,\,\cdot\,)\equiv p(\,\cdot\,,\,\omax)\equiv q(\omax,\,\cdot\,)\equiv q(\,\cdot\,,\,\omin)\equiv0
\]
whenever either $\omin$ or $\omax$ is finite. Furthermore, we assume that either $p$ and $q$ are non-zero in all other cases, or one of them is the identically zero function (totally asymmetric case).
\item In order to provide a smoothening effect in the dynamics we assume monotonicity in the following way:
\be
\ba
p(z+1,\,y)&\ge p(z,\,y),\qquad&p(y,\,z+1)&\le p(y,\,z)\\
q(z+1,\,y)&\le q(z,\,y),\qquad&q(y,\,z+1)&\ge q(y,\,z)
\ea\label{eq:mon}
\ee
for $y,\,z,\,z+1\in I$. This property has the natural interpretation that the higher neighbors a column has, the faster it grows and the slower it gets a brick removed. Our model is hence {\sl attractive}.
\item We are going to use the product property of the model's translation-invariant stationary measure. For this reason, similarly to Cocozza-Thi\-vent \cite{coco}, we need two assumptions:
\begin{itemize}
\item For any $x,\,y,\,z\in I$
\be
\ba
&p(x,\,y)+p(y,\,z)+p(z,\,x)\!\!\!\!\!&&\\
+\,&q(x,\,y)+q(y,\,z)+q(z,\,x)\!\!\!\!\!&=&\,p(x,\,z)+p(z,\,y)+p(y,\,x)\\
&&+&\,q(x,\,z)+q(z,\,y)+q(y,\,x).
\ea\label{eq:stacifelt}
\ee
\item There are symmetric functions $s_p$ and $s_q$, and a common function $f$, such that $f(\omin)=0$ whenever $\omin$ is finite, and for any $y,\,z\in I$
\be
p(y,\,z)=s_p(y,\,z+1)\cdot f(y)\qquad\text{and}\qquad q(y,\,z)=s_q(y+1,\,z)\cdot f(z).\label{eq:symm}
\ee
Condition \eqref{eq:mon} implies that $f$ is non-decreasing on $I$.
\end{itemize}
\item In order to properly construct the dynamics, restrictive growth conditions might be necessary on the rates $p$ and $q$ in case of an unbounded single-site state space $I$. We comment on this below.
\end{itemize}
At time $t$, the interface mentioned above is described by $\un{\om}(t)$. Let $\vp\,:\,\Omega\to\Rb$ be a finite cylinder function, i.e.,\ $\vp$ depends on a finite number of $\om_i$ values. The growth of this interface is a Markov process, with the formal infinitesimal generator $L$:
\be
\ba
(L\vp)(\un\om)&=\sum_{i\in\Zb}p(\om_i,\,\om_{i+1})\cdot\left[\vp(\un\om^{i,i+1})-\vp(\un\om)\right]\\
&+\sum_{i\in\Zb}q(\om_i,\,\om_{i+1})\cdot\left[\vp(\un\om^{i+1,i})-\vp(\un\om)\right].
\ea\label{eq:gen}
\ee
The construction of dynamics is available in the following situations. Several models with bounded rates are well understood and can be handled via the Hille-Yosida Theorem, see Liggett \cite{ips}. When the rates $p$ and $q$ grow at most linearly fast as functions of the local $\om$ values, then methods initiated by Liggett and Andjel lead to the construction of some zero range type systems (Andjel \cite{and}, Liggett \cite{lize},  Booth and Quant \cite{lna,lorna}). The totally asymmetric zero range and bricklayers' processes with at most exponentially growing rates are constructed in Bal\'azs, Rassoul-Agha, Sepp\"al\"ainen and Sethuraman \cite{exists}. See the definition of zero range and bricklayers' processes below.

We assume that the existence of dynamics can be established on a set of tempered configurations $\wt\Omega$ (i.e.\ configurations obeying some restrictive growth conditions), and we have the usual properties of the semigroup and the generator acting on nice functions on this set. We also assume that $\wt\Omega$ is of full measure w.r.t.\ the stationary measures defined in Section \ref{sc:gibbs}. Questions of existence of dynamics are not considered in the present paper.

\subsection{Examples}\label{sc:exa}

Many well-known nearest neighbor processes belong to this class, see e.g.\ \cite{varj2nd} for a more complete treatment. Here we only list those which we consider later in more detail.

\begin{itemize}
\item {\bf The asymmetric simple exclusion process (ASEP)} introduced by F.\ Spitzer \cite{spi} is characterised by $\omin=0,\ \omax=1$, $f(z)={\bf1}\{z=1\}$,
\[
s_p(y,\,z)=p\cdot{\bf1}\{y=z=1\}\qquad\text{and}\qquad s_q(y,\,z)=q\cdot{\bf1}\{y=z=1\},
\]
where $p>q$ are non-negative reals adding up to 1 (see \eqref{eq:symm}). In this case
\[
p(y,\,z)=p\cdot{\bf1}\{y=1,\,z=0\}\qquad\text{and}\qquad q(y,\,z)=q\cdot{\bf1}\{y=0,\,z=1\}.
\]
Here $\om_i\in\{0,\,1\}$ is the occupation number for site $i$, $p(\om_i,\,\om_{i+1})$ is the rate for a particle to jump from site $i$ to $i+1$, and $q(\om_i,\,\om_{i+1})$ is the rate for a particle to jump from site $i+1$ to $i$. These rates have values $p$ and $q$, respectively, whenever there is a particle to perform the above jumps, and there is no particle on the terminal site of the jumps. Conditions \eqref{eq:mon} and \eqref{eq:stacifelt} are also satisfied by these rates.
\item {\bf Totally asymmetric zero range processes} are included by an arbitrary nondecreasing function $f\,:\,\Zb\to\Rb^+$,
\[
\ba
s_p(y,\,z)&\equiv 1\qquad&&\text{and}&\qquad s_q(y,\,z)&\equiv0,\\
p(y,\,z)&=f(y)\qquad&&\text{and}&\qquad q(y,\,z)&\equiv0.
\ea
\]
In its original form, the totally asymmetric zero range process is a particle system with \(\om_i\) particles at site \(i\), a particle jumps from $i$ to \(i+1\) with rate $f(\om_i)$. The setting $\omin=0,\ \omax=\infty$, and so $f(0)=0$ would correspond to this situation. In the sequel it will be important for us to allow for negative values of \(\om\) as well, which comes naturally in the surface representation. We shall refer to this class of models with \(\omin=-\infty,\ \omax=\infty\), and \(f(y)>0\) for all \(y\in\Zb\), as generalized zero range processes. Conditions \eqref{eq:mon} and \eqref{eq:stacifelt} trivially hold for the rates.
\begin{itemize}
\item As a special case, the {\bf generalized totally asymmetric exponential zero range process} (we will simply refer to it  as {\bf GZRP}) is obtained by \(p(y,\,z)=f(y)=\e{\beta(y-1/2)}\) with a \(\beta>0\) parameter. Omitting the constant \(-\beta/2\) would simply correspond to a change of timescale but would bring in some unwanted factors in our final result.
\end{itemize}
\item {\bf Totally asymmetric bricklayers models.} Let \(\omin=-\infty\), \(\omax=\infty\), $f\,:\,\Zb\to\Rb^+$ non-decreasing, also having the property
\[
f(z)\cdot f(1-z)=1\qquad\text{for all }z\in\Zb.
\]
The values of $f$ for positive $z$'s thus determine the values for non-positive $z$'s. Set
\[
s_p(y,\,z)=1+\frac{1}{f(y)f(z)},\qquad s_q(y,\,z)\equiv0,
\]
which results in
\[
p(y,\,z)=f(y)+f(-z),\qquad q(y,\,z)\equiv0.
\]
This process was first represented by bricklayers standing at each site $i$, laying a brick on the column on their left with rate $f(-\om_i)$ and laying a brick to their right with rate $f(\om_i)$, hence the name. Conditions \eqref{eq:mon} and \eqref{eq:stacifelt} hold for the rates.
\begin{itemize}
\item As a special case, the {\bf totally asymmetric exponential bricklayers process} (we will just abbreviate as {\bf BLP}) is obtained by \(f(z)=\e{\beta(z-1/2)}\) with a \(\beta>0\) parameter.
\end{itemize}
Note that the BLP is a symmetrized version of the GZRP in the following sense. For a GZRP with a given $\beta$ one can consider its counterpart by a space reflection in the surface representation. The symmetric combination of these two processes (which is obtained by taking the sum of the two generators) gives the corresponding BLP (in the particle representation a particle-antiparticle transformation should follow the space reflection).
\end{itemize}

\subsection{Translation invariant stationary product distributions}\label{sc:gibbs}

We now present some translation invariant stationary distributions for these processes. For many cases it has been proved that these are the only extremal translation-invariant stationary distributions. Following some ideas in Cocozza-Thivent \cite{coco}, we first consider the non-decreasing function $f$ of \eqref{eq:symm}.
 For $I\ni z>0$ we define 
\[
f(z)!:\,=\prod_{y=1}^zf(y),
\]
while for $I\ni z<0$ let
\[
f(z)!:\,=\frac{1}{\prod\limits_{y=z+1}^0f(y)},
\]
finally $f(0)!:\,=1$. Then we have 
\[
f(z)!\cdot f(z+1)=f(z+1)!
\]
for all $z\in I$. Let
\[
\bar\te:\,=\left\{\begin{array}{ll}\lim\limits_{z\to\infty}\log(f(z))\ \ &,\ \text{if}\ \omax=\infty\\\infty\ \ &,\ \text{else}\end{array}\right.
\]
and
\[
\un\te:\,=\left\{\begin{array}{ll}\lim\limits_{z\to\infty}\log(f(-z))\ \ &,\ \text{if}\ \omin=-\infty\\-\infty\ \ &,\ \text{else}.\end{array}\right.
\]
By monotonicity of $f$, we have $\bar\te\ge\un\te$. We assume $\bar\te>\un\te$. With a generic real parameter $\te\in\left(\un\te,\,\bar\te\right)$, which is often referred to as the chemical potential, we define the partition sum as
\[
Z(\te):\,=\sum_{z\in I}\frac{\e{\te z}}{f(z)!}<\infty.
\]
Let the product-distribution $\un\mu^\te$ have marginals
\be
\mu^\te(z)=\un\mu^\te\left\{\un\om\,:\,\om_i=z\right\}:\,=\left\{
\ba
&\frac{1}{Z(\te)}\cdot\frac{\e{\te z}}{f(z)!}&&\text{if }z\in I,\\
&0&&\text{if }z\notin I.
\ea
\right.\label{eq:om}
\ee
Then the product distribution $\un\mu^\te$ is stationary for the process generated by \eqref{eq:gen}.

We define the function
\be
\rho(\te):\,=\frac{\di}{\di\te}\log(Z(\te))=\sum_{z\in I}z\cdot\mu^\te(z),\qquad\un\te<\te<\bar\te.\label{eq:rote}
\ee
This is the density of particles, and it is quite easy to see that this is a strictly increasing function. Its inverse will be denoted by \(\te(\vr)\). Due to the bijection it is always possible to choose a suitable characterization of the stationary state either with the chemical potential $\theta$ or with the corresponding density \(\vr\).

As for our examples, the ASEP has a Bernoulli product stationary distribution of which the density parameter is
\be
\rho(\te)=\frac{\e{\te}}{1+\e{\te}}.\label{eq:asepro}
\ee

In the case of both our exponential GZRP and BLP examples, computing the factorials in \eqref{eq:om} results in the discrete Gaussian
\[
\mu^\te(z)=\frac{\e{\te^2\!/2\beta}}{Z(\te)}\cdot\e{-\frac\beta2\cdot\bigl(z-\frac\te\beta\bigr)^2}\qquad z\in\Zb,
\]
from which
\[
Z(\te)=\e{\te^2\!/2\beta}\cdot\sum_{z=-\infty}^\infty\e{-\frac\beta2\cdot\bigl(z-\frac\te\beta\bigr)^2}.
\]
While no explicit form is available for the partition sum $Z(\te)$, the identity
\be
Z(\te-\beta)=\e{\beta/2-\te}\cdot Z(\te)\label{eq:expz}
\ee
can easily be shown. Via \eqref{eq:rote} this implies
\be
\rho(\te-\beta)=\rho(\te)-1\label{eq:expkul}
\ee
for the exponential GZRP and BLP processes. Notice also that
\be
\mu^{\te-\beta}(y)=\mu^\te(y+1).\label{eq:mutol}
\ee

\subsection{Hydrodynamics, very briefly}

In preparation for stating the result, we first briefly mention that it is believed, and in many cases proved, that models of our family satisfy a conservation law of the form
\be
\partial_T\rho(T,\,X)+\partial_X\Hc(\rho(T,\,X))=0\label{eq:hcl}
\ee
in the \emph{Eulerian scaling}, with the density \(\rho\) being a function of the rescaled time and space variables \(T\) and \(X\), and \(\Hc(\vr)\), the flux function, being the expected net current in the  stationary measure with density \(\vr\) (recall \eqref{eq:rote}):
\[
\Hc(\vr)=\sum_{y,z\in I}[p(y,\,z)-q(y,\,z)]\cdot\mu^{\te(\vr)}(y)\mu^{\te(\vr)}(z).
\]
See e.g.\ Rezakhanlou \cite{hl} or Bahadoran, Guiol, Ravishankar and Saada \cite{bagurasa} for details.

At time zero let the density profile be a step function: \(\rho(0,\,X)\equiv\vr\) for \(X<0\) and \(\rho(0,\,X)\equiv\la\) for \(X>0\). This initial condition corresponds to a \emph{shock} if \(\vr<\la\) and \(\Hc\) is concave, or if \(\vr>\la\) and \(\Hc\) is convex. The entropy solutions of \eqref{eq:hcl} in these cases are rigid translations of the shock with the \emph{Rankine-Hugoniot velocity}
\be
V=\frac{\Hc(\la)-\Hc(\vr)}{\la-\vr}.\label{eq:rh}
\ee 
From this formula one can check that multiple shocks eventually meet and merge into a single shock having the leftmost and the rightmost of the initial density values on its left and right sides, respectively.

It is easy to see that for the ASEP we have
\[
\Hc(\vr)=(p-q)\cdot\vr(1-\vr)
\]
which is a concave function (\(p>q\)). It is also straightforward from the definitions that
\[
\Hc(\vr)=\e{\te(\vr)}\qquad\text{and}\qquad\Hc(\vr)=\e{\te(\vr)}+\e{-\te(\vr)}
\]
for the GZRP and BLP examples, respectively. Due to the lack of an explicit formula for \(\te(\vr)\), it is a nontrivial fact that both these flux functions are convex when the rate functions \(f\) are convex \cite{convex}.

\subsection{The second class particle}

Attractivity \eqref{eq:mon} makes it possible to define the \emph{second class particle}. Let \(\un\de^j\) be the vector of components 1 for site \(j\), and 0 for all other sites in \(\Zb\). Let \(\un\om\in\wt\Omega\) such that \(\om_j<\omax\), and
\be
\un\ze:\,=\un\om+\un\de^j.\label{eq:deltakul}
\ee
We say that there is a second class particle at site \(j\). We let the pair \((\un\om,\,\un\ze)\) evolve in the \emph{basic coupling}%
. That is, assuming and making use of \eqref{eq:deltakul}, the effect of the coupled generator (which we also denote by \(L\)) for a finite cylinder function \(\vp\) of a pair \((\un\om,\,\un\ze)\) is
\begin{align}
(L\vp)(\un\om,\,\un\ze)&=\sum_{i=a-1}^{j-2}p(\om_i,\,\om_{i+1})\cdot\bigl[\vp(\un\om^{i,i+1},\,\un\ze^{i,i+1})-\vp(\un\om,\,\un\ze)\bigr]\label{eq:cougen}\\
&+\sum_{i=j+1}^bp(\om_i,\,\om_{i+1})\cdot\bigl[\vp(\un\om^{i,i+1},\,\un\ze^{i,i+1})-\vp(\un\om,\,\un\ze)\bigr]\notag\\
&\quad+p(\om_{j-1},\,\om_j+1)\cdot\bigl[\vp(\un\om^{j-1,j},\,\un\ze^{j-1,j})-\vp(\un\om,\,\un\ze)\bigr]\notag\\
&\quad+\bigl[p(\om_{j-1},\,\om_j)-p(\om_{j-1},\,\om_j+1)\bigr]\cdot\bigl[\vp(\un\om^{j-1,j},\,\un\ze)-\vp(\un\om,\,\un\ze)\bigr]\label{eq:ombal}\\
&\quad+p(\om_j,\,\om_{j+1})\cdot\bigl[\vp(\un\om^{j,j+1},\,\un\ze^{j,j+1})-\vp(\un\om,\,\un\ze)\bigr]\notag\\
&\quad+\bigl[p(\om_j+1,\,\om_{j+1})-p(\om_j,\,\om_{j+1})\bigr]\cdot\bigl[\vp(\un\om,\,\un\ze^{j,j+1})-\vp(\un\om,\,\un\ze)\bigr]\label{eq:zejobb}\\
&+\sum_{i=a-1}^{j-2}q(\om_i,\,\om_{i+1})\cdot\bigl[\vp(\un\om^{i+1,i},\,\un\ze^{i+1,i})-\vp(\un\om,\,\un\ze)\bigr]\notag\\
&+\sum_{i=j+1}^bq(\om_i,\,\om_{i+1})\cdot\bigl[\vp(\un\om^{i+1,i},\,\un\ze^{i+1,i})-\vp(\un\om,\,\un\ze)\bigr]\notag\\
&\quad+q(\om_{j-1},\,\om_j)\cdot\bigl[\vp(\un\om^{j,j-1},\,\un\ze^{j,j-1})-\vp(\un\om,\,\un\ze)\bigr]\notag\\
&\quad+\bigl[q(\om_{j-1},\,\om_j+1)-q(\om_{j-1},\,\om_j)\bigr]\cdot\bigl[\vp(\un\om,\,\un\ze^{j,j-1})-\vp(\un\om,\,\un\ze)\bigr]\label{eq:zebal}\\
&\quad+q(\om_j+1,\,\om_{j+1})\cdot\bigl[\vp(\un\om^{j+1,j},\,\un\ze^{j+1,j})-\vp(\un\om,\,\un\ze)\bigr]\notag\\
&\quad+\bigl[q(\om_j,\,\om_{j+1})-q(\om_j+1,\,\om_{j+1})\bigr]\cdot\bigl[\vp(\un\om^{j+1,j},\,\un\ze)-\vp(\un\om,\,\un\ze)\bigr]\label{eq:omjobb},
\end{align}
where \(\vp\) depends on the configuration over the sites \(a\dots b\). We assume that \(a<-1<1<b\). This generator gives the correct marginal evolution for each \(\un\om\) and \(\un\ze\). Only steps \eqref{eq:ombal} - \eqref{eq:omjobb} influence the second class particle. In fact these steps result in jumps of this particle, hence the single second class particle is conserved for all times.

Constructing the dynamics of a coupled pair with any number of second class particles can be done along the same lines, and their number is again conserved, see the generator for that case in Section \ref{sc:mul}.

\section{Shock measures with second class particles}\label{sc:results}

We now define product shock measures with a single second class particle on a coupled pair of processes. Later on multiple shocks and second class particles will be considered, but the case of one second class particle is much simpler, so we first demonstrate it here. The marginals will be like the \eqref{eq:om} stationary ones, except for one site which has the second class particle. Define, for \(\un\te<\te<\bar\te\), the one-site marginal \(\nu^\te\) on \(I\times I\) by
\[
\nu^\te(y,\,z)=\left\{\ba
&\mu^\te(y),&&\text{if }y=z,\\
&0,&&\text{if }y\ne z,
\ea\right.
\]
where \(\mu^\te\) is the stationary marginal \eqref{eq:om}. Let also \(\hat\mu\) be a measure on \(I\) such that \(\hat\mu(\omax)=0\) if \(\omax\) is finite. Define
\[
\hat\nu(y,\,z)=\left\{\ba
&\hat\mu(y),&&\text{if }z=y+1,\\
&0,&&\text{otherwise.}
\ea\right.
\]
With these marginals and with \(\un\te<\te,\,\si<\bar\te\) we define the product measure
\be
\un\nu_j:\,=\bigotimes_{i<j}\nu^\te\bigotimes_{i=j}\hat\nu\bigotimes_{i>j}\nu^\si\label{eq:shm}
\ee
of marginals \(\nu^\te\) on the left of site \(j\), \(\hat\nu\) at \(j\), and \(\nu^\si\) on the right of site \(j\). Therefore, \(\un\nu_j\) is a measure on coupled pairs with exactly one second class particle at site \(j\). 

In the case of the ASEP it is more convenient to use the densities as parameters instead of fugacities. Here, with a slight abuse of notation, the results are expressed in terms of $\vr=\rho(\theta)$ and $\lambda=\rho(\sigma)$ being the left and right densities correspondingly, according to \eqref{eq:asepro}.

\subsection{Results for a single shock}

The main result is on the time-evolution of some particular distributions of a coupled pair \((\un\om(t),\,\un\ze(t))\) of some particular models. We use the semigroup notation for the evolution of such a distribution, so \(\un\nu S(t)\) denotes the distribution of the pair at time \(t\) if the initial distribution at time \(0\) was \(\un\nu\). Recall the definition \eqref{eq:shm}.
\begin{tm}\label{tm:main}
The identity
\be
\frac{\di}{\di t}\un\nu_jS(t)\Bigr|_{t=0}=P\cdot[\un\nu_{j+1}-\un\nu_j]+Q\cdot[\un\nu_{j-1}-\un\nu_j]\label{eq:main}
\ee
holds, in the sense of test functions, in the following special cases with the following parameters \(P\) and \(Q\):
\begin{itemize}
\item For the ASEP, if the relation
\be
\frac{\la(1-\vr)}{\vr(1-\la)}=\frac pq\label{eq:aseprel}
\ee
holds between the densities and the asymmetry, and
\be
\hat\mu(0)=1.\label{eq:asepmu}
\ee
In this case
\be
\ba
P&=\frac{1-\la}{1-\vr}\cdot p=\frac\la\vr\cdot q=(1-\la)p+\la q\qquad\text{and}\\
Q&=\frac{1-\vr}{1-\la}\cdot q=\frac\vr\la\cdot p=(1-\vr)q+\vr p.
\ea\label{eq:aseppq}
\ee
\item For the exponential GZRP and BLP as defined in Section \ref{sc:exa}, if the relation
\be
\te-\si=\beta\label{eq:exprel}
\ee
holds between the parameters, and
\be
\hat\mu(y)=\mu^\si(y).\label{eq:expmu}
\ee
In this case we have
\be
\ba
P&=\e{\te}-\e{\si}\qquad\text{and}&\qquad Q&=0&\qquad&\text{for the GZRP},\\
P&=\e{\te}-\e{\si}\qquad\text{and}&\qquad Q&=\e{-\si}-\e{-\te}&\qquad&\text{for the BLP}.
\ea\label{eq:blppq}
\ee
\end{itemize}
\end{tm}
The only way to have a second class particle at \(j\) in ASEP is \(\om_j=0\) and \(\ze_j=1\). This fact is of course reflected in the form of \(\hat\nu\). This is not the case for GZRP and BLP, where the marginals for \(\om_j\) and \(\ze_j\) are not a priori restricted in any way. The result turns out to require these marginals to be \(\mu^\si\) and \(\mu^\te\), respectively (this latter being a consequence of \eqref{eq:mutol}).
\begin{rem}
As in \cite{qse} and \cite{sokvalak}, \eqref{eq:main} has a natural random walk interpretation: the shock measure \(\un\nu_j\), with the second class particle in the middle, performs an asymmetric simple random walk. The drift \(P-Q\) of this walk agrees with the Rankine Hugoniot velocity \eqref{eq:rh}. This can be seen for the ASEP from the definitions and using the relation \eqref{eq:aseprel}. For the exponential GZRP and BLP models the relation \eqref{eq:exprel} together with \eqref{eq:expkul} can be used to check this.
\end{rem}

We remark here that a theorem, similar in spirit, will be shown in Section \ref{sc:bcrwres} for the branching coalescing random walk model.

In \cite{qse} a shock measure with density \(\vr\) for sites at or left of \(j\) and \(\la\) on the right of site \(j\) was first considered. It was shown there that this structure performs a random walk with the above jump rates. In view of our result, mixing the two marginals of \(\un\nu_j\) leads to \cite{qse}'s shock measure: flip a biased coin initially, and with probability \(\la\) follow the distribution of the upper process \(\un\ze\) while with probability \(1-\la\) follow the lower process \(\om\). Later on, in their Theorem 2, Belitsky and Sch\"utz in \cite{qse} put particles at the positions of the shocks to describe interaction (see our Section \ref{sc:mul} for interaction of shocks), and they also remark that these are \emph{not} the second class particles. However, their setting is also obtained from our results by considering the second marginal \(\un\ze\) of \(\un\nu_j\). Therefore, their particles at the shock positions can also be viewed as second class particles. Later on, the idea of random walking shocks \emph{with} second class particles at the shock positions appeared in the form of a conjecture in \cite{bcrw} in the context of the ASEP with open boundaries.

We recover the results of \cite{sokvalak} if we only consider the first marginal \(\un\om\) of \(\un\nu_j\) for the exponential BLP.

\subsection{Results for multiple shocks}\label{sc:mul}

While we tried to keep notations as simple as possible for a single shock, multiple shocks seem to require a more complicated treatment. We describe below and use a formalism capable of handling multiple shocks with second class particles. First we rewrite \eqref{eq:cougen} in the general case, only assuming \(\om_i\le\ze_i\) for each \(i\in\Zb\):
\begin{align}
(L\vp)(\un\om,\,\un\ze)=\sum_{i=a-1}^b\Bigl\{&p(\om_i,\,\ze_{i+1})\cdot\bigl[\vp(\un\om^{i,i+1},\,\un\ze^{i,i+1})-\vp(\un\om,\,\un\ze)\bigr]\label{eq:altcougen}\\
+\bigl[&p(\om_i,\,\om_{i+1})-p(\om_i,\,\ze_{i+1})\bigr]\cdot\bigl[\vp(\un\om^{i,i+1},\,\un\ze)-\vp(\un\om,\,\un\ze)\bigr]\label{eq:altomb}\\
+\bigl[&p(\ze_i,\,\ze_{i+1})-p(\om_i,\,\ze_{i+1})\bigr]\cdot\bigl[\vp(\un\om,\,\un\ze^{i,i+1})-\vp(\un\om,\,\un\ze)\bigr]\label{eq:altzej}\\
+\phantom{\bigl[}&q(\ze_i,\,\om_{i+1})\cdot\bigl[\vp(\un\om^{i+1,i},\,\un\ze^{i+1,i})-\vp(\un\om,\,\un\ze)\bigr]\notag\\
+\bigl[&q(\ze_i,\,\ze_{i+1})-q(\ze_i,\,\om_{i+1})\bigr]\cdot\bigl[\vp(\un\om,\,\un\ze^{i+1,i})-\vp(\un\om,\,\un\ze)\bigr]\Bigr\}\label{eq:altzeb}\\
+\bigl[&q(\om_i,\,\om_{i+1})-q(\ze_i,\,\om_{i+1})\bigr]\cdot\bigl[\vp(\un\om^{i+1,i},\,\un\ze)-\vp(\un\om,\,\un\ze)\bigr]\label{eq:altomj}.
\end{align}
Again, the above differences of jump rates are nonnegative, and marginally both \(\un\om\) and \(\un\ze\) evolve according to the original dynamics. Moreover, the generator keeps \(\om_i\le\ze_i\) for each site.

Let \(m\ge0\) be an integer and \(\un\te<\si<\bar\te\) a parameter value. Define the one-site marginal of a coupled pair \(\un\om,\,\un\ze\) by
\[
\nu^{\si,m}(y,\,z)=\left\{\ba
&\hat\mu^{\si,m}(y),&&\text{if }z=y+m,\\
&0,&&\text{if }z\ne y+m.
\ea\right.
\]
In our examples this marginal will coincide with the stationary marginal \eqref{eq:om} if \(m=0\). When \(m>0\), it describes a distribution on a site with \(m\) second class particles; \(\om=\ze-m\) then has distribution \(\hat\mu^{\si,m}\). We require that \(\hat\mu^{\si,m}\) gives probability zero on values \(y>\omax-m\) when \(\omax\) is finite. Now with a vector \(\un m\) of nonnegative integer components \(m_i\) and a parameter vector \(\un\si\), define the product measure
\be
\un\nu^{\un\si,\un m}:\,=\bigotimes_{i\in\Zb}\nu^{\si_i,m_i}\label{eq:unnudef}
\ee
on pairs \((\un\om,\,\un\ze)\). This measure describes \(m_i\) second class particles at site \(i\) with marginal \(\hat\mu^{\si_i,m_i}\) for \(\om_i=\ze_i-m_i\) and, in our examples, stationary marginals \(\hat\mu^{\si_i,0}=\mu^{\si_i}\) at other sites with no second class particles. We assume a finite number of second class particles in the system that is, \(\sum\limits_{i\in\Zb}m_i<\infty\).

In the ASEP case we still use \eqref{eq:asepro} and switch to the \(\vr_i\) parametrization. We also adopt the notation \eqref{eq:ugrdef} for \(\un m^{i,j}\), while \(\un\si^{i,j}\) will be some modified parameter vector which we specify below in the theorem:
\begin{tm}\label{tm:mush}
The identity
\be
\ba
\frac{\di}{\di t}\un\nu^{\un\si,\un m}S(t)\Bigr|_{t=0}=\sum_{i\in\Zb}&P(m_i,\,m_{i+1},\,\si_i,\,\si_{i+1})\cdot\bigl[\un\nu^{\un\si^{i,i+1},\un m^{i,i+1}}-\un\nu^{\un\si,\un m}\bigr]\\
+\,&Q(m_i,\,m_{i+1},\,\si_i,\,\si_{i+1})\cdot\bigl[\un\nu^{\un\si^{i+1,i},\un m^{i+1,i}}-\un\nu^{\un\si,\un m}\bigr]
\ea\label{eq:mush}
\ee
holds, in the sense of test functions, in the following special cases with the following parameters \(\un\si^{i,i+1}\), \(\un\si^{i+1,i}\), \(P(m_i,\,m_{i+1},\,\si_i,\,\si_{i+1})\) and \(Q(m_i,\,m_{i+1},\,\si_i,\,\si_{i+1})\):
\begin{itemize}
\item For the ASEP, if \(m_i=0\) or \(1\), and the relation
\be
\frac{\vr_{i+1}(1-\vr_i)}{\vr_i(1-\vr_{i+1})}=\left\{\ba
&1,&&\text{if }m_{i+1}=0,\\
&\frac pq,&&\text{if }m_{i+1}=1
\ea\right.\label{eq:asepsokfelt}
\ee
holds between the densities and the asymmetry, and
\be
\hat\mu^{\vr,m}(0)=1-\hat\mu^{\vr,m}(1)=\left\{\ba
&1-\vr,&&\text{if }m=0,\\
&1,&&\text{if }m=1.
\ea\right.\label{eq:asepsokmuhat}
\ee
In this case
\begin{align}
\vr_j^{i,i+1}&=\left\{\ba
&\vr_j,&&\text{for }j\ne i,\\
&\vr_{i-1},&&\text{for }j=i,
\ea\right.\label{eq:asepvrj}\\
\vr_j^{i+1,i}&=\left\{\ba
&\vr_j,&&\text{for }j\ne i,\\
&\vr_{i+1},&&\text{for }j=i\text{, and}
\ea\right.\label{eq:asepvrb}\\
P(m_i,\,m_{i+1},\,\vr_i,\,\vr_{i+1})&=m_i(1-m_{i+1})\cdot\bigl[(1-\vr_{i+1})p+\vr_{i+1}q\bigr],\label{eq:asepshrj}\\
Q(m_i,\,m_{i+1},\,\vr_i,\,\vr_{i+1})&=(1-m_i)m_{i+1}\cdot\bigl[(1-\vr_i)q+\vr_ip\bigr].\label{eq:asepshrb}
\end{align}
\item For the exponential GZRP and BLP as defined in Section \ref{sc:exa}, if the relation
\be
\si_{i-1}-\si_i=\beta m_i\label{eq:blpsokfelt}
\ee
holds between the parameters, and
\[
\hat\mu^{\si,m}(y)=\mu^\si(y),
\]
the stationary marginal \eqref{eq:om}, regardless of \(m\). In this case we have
\begin{align}
\si_j^{i,i+1}&=\left\{\ba
&\si_j,&&\text{for }j\ne i,\\
&\si_i+\beta,&&\text{for }j=i,
\ea\right.\label{eq:blpsij}\\
\si_j^{i+1,i}&=\left\{\ba
&\si_j,&&\text{for }j\ne i,\\
&\si_i-\beta,&&\text{for }j=i,
\ea\right.\label{eq:blpsib}\\
P(m_i,\,m_{i+1},\,\si_i,\,\si_{i+1})&=\e{\si_i+\beta m_i}-\e{\si_i}&&\text{and}\label{eq:zrpshrj}\\
Q(m_i,\,m_{i+1},\,\si_i,\,\si_{i+1})&=0&&\quad\text{for the GZRP},\\
P(m_i,\,m_{i+1},\,\si_i,\,\si_{i+1})&=\e{\si_i+\beta m_i}-\e{\si_i}&&\text{and}\\
Q(m_i,\,m_{i+1},\,\si_i,\,\si_{i+1})&=\e{-\si_{i+1}}-\e{-\si_{i+1}-\beta m_{i+1}}&&\quad\text{for the BLP}.\label{eq:blpshrb}
\end{align}
\end{itemize}
\end{tm}
First notice that \eqref{eq:asepvrj} - \eqref{eq:asepvrb} and \eqref{eq:blpsij} - \eqref{eq:blpsib} keep \eqref{eq:asepsokfelt} and \eqref{eq:blpsokfelt}, respectively, valid after each jump. Also notice that this theorem reduces to Theorem \ref{tm:main} in the case of a single shock. In \eqref{eq:zrpshrj} - \eqref{eq:blpshrb}, the exponents \(\si_j+\beta m_j\) are simply the parameter values of the \(\ze_j\) marginal or the coupled pair \((\om_j,\,\ze_j)\) due to \eqref{eq:mutol}.

In both the ASEP and the GZRP/BLP, the result can be given a simple interpretation. Thinking about the process \(\un m(t)\) of second class particles, the right-hand side \eqref{eq:mush} is a generator of the interacting system \(\un m(t)\). Then the evolution rules for the densities keep the necessary conditions \eqref{eq:asepsokfelt} and \eqref{eq:blpsokfelt} valid by adjusting the density values to \(\un m(t)\).

For the ASEP \eqref{eq:mush}, together with \eqref{eq:asepshrj} - \eqref{eq:asepshrb}, looks like an asymmetric exclusion process except that the (annealed) jump rates of the second class particles depend on the density values. Due to the exclusion rule, these particles always keep their order, and \eqref{eq:asepshrj} - \eqref{eq:asepshrb} simply say that their right jump rates decrease and left jump rates increase as we go from the leftmost to the rightmost second class particle. (Recall that \(p>q\) and the density increases from left to right at each second class particle.) The second class particles thus stay within a tight distance from each other, which is already invisible on the macroscopic scale of the conservation law. We show in the Appendix that the velocity of this structure of shocks agrees with the Rankine Hugoniot velocity \eqref{eq:rh} taken with the leftmost and the rightmost density values. Once the particles, put in the shock positions in \cite{qse}, are replaced with our second class particles, our result coincides with the interaction of the random walking shocks of \cite{qse}.

For the exponential GZRP/BLP, \eqref{eq:mush} with the rates \eqref{eq:zrpshrj} - \eqref{eq:blpshrb} becomes a zero range-type generator i.e., the (annealed) jump rate of a second class particle only depends on the local configuration at its site. However, these jump rates again depend, via the parameter values \(\un\si\), on the number of other second class particles in front and behind. This dependence is such that the more behind, the greater right jump rate and smaller left jump rate for a second class particle. Lemma 4.3 from \cite{sokvalak} is word for word valid also with second class particles: the \(n\) particles stay within a tight distance from each other, and their center of mass performs a drifted simple random walk on the lattice \(\Zb/n\) with right jump rate \(\e{\si_\text{left}}-\e{\si_\text{right}}\) and left jump rate \(\e{-\si_\text{right}}-\e{-\si_\text{left}}\) for BLP and zero for GZRP. Here \(\si_\text{left}\) is the leftmost, \(\si_\text{right}\) is the rightmost value of \(\si_i\). The resulting drift of the center of mass of course coincides with the Rankine-Hugoniot velocity for a large shock of densities \(\rho(\si_\text{left})\) and \(\rho(\si_\text{right})\).

\section{Proofs}\label{sc:proofs}

\subsection{The proof for a single shock}

\begin{proof}[Proof of Theorem \ref{tm:main}]
We assume, without loss of generality, that \(j=0\). To prove Theorem \ref{tm:main} we take a finite cylinder function \(\vp\) on \(\wt\Omega^2\) and see how the coupled generator \eqref{eq:cougen} acts on \(\vp\) under \(\Ev_0\), the expectation w.r.t.\ \(\un\nu_0\). This expectation involves summations for all \(\om_i\) variables. We shall change the summation variables in each term where the argument of \(\vp\) is modified such that, after the change, each term will contain \(\vp(\un\om,\,\un\ze)\).

Lines \eqref{eq:ombal} - \eqref{eq:omjobb} describe jumps of the second class particle. For all the other lines, the second class particle stays at site \(0\), and the change of variables can be done without much complications. In the term with \(\vp(\un\om^{-1,0},\,\un\ze^{-1,0})\) in the second line, for example, we will put \(\un\om^{-1,0}\) and \(\un\ze^{-1,0}\) as our new summation variables. This step will change the arguments of the rates, and it will also bring in factors of \(\mu^\te\) and \(\hat\mu\).

For lines \eqref{eq:ombal} - \eqref{eq:omjobb}, however, there is a complication. While \(\Ev_0\) gives probability one on \(\ze_0=\om_0+1\), this will not hold after the change of variables. As an example, the new variables for the term with \(\vp(\un\om^{-1,0},\,\un\ze)\) in line \eqref{eq:ombal} are \(\un\om^{-1,0}\) and \(\un\ze\), for which we rather have
\[
\ze_{-1}=\om^{-1,0}_{-1}+1\qquad\text{and}\qquad\ze_0=\om^{-1,0}_0.
\]
This is the manifestation of the fact that \eqref{eq:ombal} describes a left jump of the second class particle from site \(0\) to site \(-1\). Notice that this new configuration \((\un\om^{-1,0},\,\un\ze)\) is now singular to the measure \(\un\nu_0\), but not to \(\un\nu_{-1}\). Therefore the change of variables will also include a change in the expectation from \(\Ev_0\) to \(\Ev_{-1}\) for that term in line \eqref{eq:ombal}.

The result of changing variables will be of the form
\be
\Ev_0(L\vp)(\un\om,\,\un\ze)=\Ev_0\bigl\{\vp(\un\om,\,\un\ze)\cdot A\bigr\}+\Ev_{-1}\bigl\{\vp(\un\om,\,\un\ze)\cdot[B+D]\bigr\}+\Ev_1\bigl\{\vp(\un\om,\,\un\ze)\cdot[C+E]\bigr\},\label{eq:ae}
\ee
where \(B\) comes from the term with \(\vp(\un\om^{-1,0},\,\un\ze)\) of line \eqref{eq:ombal}, \(C\) comes from the term with \(\vp(\un\om,\,\un\ze^{0,1})\) of line \eqref{eq:zejobb}, \(D\) comes from the term with \(\vp(\un\om,\,\un\ze^{0,-1})\) of line \eqref{eq:zebal}, \(E\) comes from the term with \(\vp(\un\om^{1,0},\,\un\ze)\) of line \eqref{eq:omjobb}, and \(A\) from all other terms. The aim is to show that \(A\), \(B+D\), and \(C+E\) each do not depend on the variables \(\om_a\dots\om_b\). In fact, \(B+D\) will give \(Q\), \(C+E\) will amount to \(P\), while \(A\) will make up for the negative terms \(-P-Q\) in \eqref{eq:main}. We now compute the quantities \(A\) - \(E\).

We start with \(A\) and write it line for line as terms come from \eqref{eq:cougen}:
\[
\ba
A&=\sum_{i=a-1}^{-2}\bigl[p(\om_i+1,\,\om_{i+1}-1)\cdot\frac{\mu^\te(\om_i+1)\mu^\te(\om_{i+1}-1)}{\mu^\te(\om_i)\mu^\te(\om_{i+1})}-p(\om_i,\,\om_{i+1})\bigr]\\
&\quad+\sum_{i=1}^b\bigl[p(\om_i+1,\,\om_{i+1}-1)\cdot\frac{\mu^\si(\om_i+1)\mu^\si(\om_{i+1}-1)}{\mu^\si(\om_i)\mu^\si(\om_{i+1})}-p(\om_i,\,\om_{i+1})\bigr]\\
&\quad+p(\om_{-1}+1,\,\om_0)\cdot\frac{\mu^\te(\om_{-1}+1)\hat\mu(\om_0-1)}{\mu^\te(\om_{-1})\hat\mu(\om_0)}-p(\om_{-1},\,\om_0+1)\\
&\quad-p(\om_{-1},\,\om_0)+p(\om_{-1},\,\om_0+1)\\
&\quad+p(\om_0+1,\,\om_1-1)\cdot\frac{\hat\mu(\om_0+1)\mu^\si(\om_1-1)}{\hat\mu(\om_0)\mu^\si(\om_1)}-p(\om_0,\,\om_1)\\
&\quad-p(\om_0+1,\,\om_1)+p(\om_0,\,\om_1)\\
&\quad+\sum_{i=a-1}^{-2}\bigl[q(\om_i-1,\,\om_{i+1}+1)\cdot\frac{\mu^\te(\om_i-1)\mu^\te(\om_{i+1}+1)}{\mu^\te(\om_i)\mu^\te(\om_{i+1})}-q(\om_i,\,\om_{i+1})\bigr]\\
&\quad+\sum_{i=1}^b\bigl[q(\om_i-1,\,\om_{i+1}+1)\cdot\frac{\mu^\si(\om_i-1)\mu^\si(\om_{i+1}+1)}{\mu^\si(\om_i)\mu^\si(\om_{i+1})}-q(\om_i,\,\om_{i+1})\bigr]\\
&\quad+q(\om_{-1}-1,\,\om_0+1)\cdot\frac{\mu^\te(\om_{-1}-1)\hat\mu(\om_0+1)}{\mu^\te(\om_{-1})\hat\mu(\om_0)}-q(\om_{-1},\,\om_0)\\
&\quad-q(\om_{-1},\,\om_0+1)+q(\om_{-1},\,\om_0)\\
&\quad+q(\om_0,\,\om_1+1)\cdot\frac{\hat\mu(\om_0-1)\mu^\si(\om_1+1)}{\hat\mu(\om_0)\mu^\si(\om_1)}-q(\om_0+1,\,\om_1)\\
&\quad-q(\om_0,\,\om_1)+q(\om_0+1,\,\om_1).
\ea
\]
Then
\be
\ba
B+D&=\bigl[p(\om_{-1}+1,\,\om_0-1)-p(\om_{-1}+1,\,\om_0)\bigr]\cdot\frac{\mu^\te(\om_{-1}+1)\hat\mu(\om_0-1)}{\hat\mu(\om_{-1})\mu^\si(\om_0)}\\
&\quad+\bigl[q(\om_{-1},\,\om_0+1)-q(\om_{-1},\,\om_0)\bigr]\cdot\frac{\mu^\te(\om_{-1})\hat\mu(\om_0)}{\hat\mu(\om_{-1})\mu^\si(\om_0)},
\ea\label{eq:bd}
\ee
and
\be
\ba
C+E&=\bigl[p(\om_0+1,\,\om_1)-p(\om_0,\,\om_1)]\cdot\frac{\hat\mu(\om_0)\mu^\si(\om_1)}{\mu^\te(\om_0)\hat\mu(\om_1)}\\
&\quad+\bigl[q(\om_0-1,\,\om_1+1)-q(\om_0,\,\om_1+1)\bigr]\cdot\frac{\hat\mu(\om_0-1)\mu^\si(\om_1+1)}{\mu^\te(\om_0)\hat\mu(\om_1)}.
\ea\label{eq:ce}
\ee
We now simplify the expression of \(A\) by plugging in the fraction of \(\mu^\te\)'s from \eqref{eq:om}. In the summations we also make use of
\[
\ba
p(\om_i+1,\,\om_{i+1}-1)\cdot\frac{f(\om_{i+1})}{f(\om_i+1)}&=p(\om_{i+1},\,\om_i)\qquad\text{and}\\
q(\om_i-1,\,\om_{i+1}+1)\cdot\frac{f(\om_i)}{f(\om_{i+1}+1)}&=q(\om_{i+1},\,\om_i),
\ea
\]
a consequence of \eqref{eq:symm}, and then
\begin{multline}
\ba
&p(\om_i,\,\om_{i-1})-p(\om_{i-1},\,\om_i)+p(\om_{i+1},\,\om_i)-p(\om_i,\,\om_{i+1})\\
+\,&q(\om_i,\,\om_{i-1})-q(\om_{i-1},\,\om_i)+q(\om_{i+1},\,\om_i)-q(\om_i,\,\om_{i+1})
\ea
\\
=p(\om_{i+1},\,\om_{i-1})-p(\om_{i-1},\,\om_{i+1})+q(\om_{i+1},\,\om_{i-1})-q(\om_{i-1},\,\om_{i+1}),\label{eq:telesc}
\end{multline}
which is a rewriting of \eqref{eq:stacifelt}. This latter produces a telescopic sum.
% \[
% \ba
% A&=p(\om_{-1},\,\om_{a-1})-p(\om_{a-1},\,\om_{-1})\\
% &\quad+p(\om_{b+1},\,\om_1)-p(\om_1,\,\om_{b+1})\\
% &\quad+p(\om_{-1}+1,\,\om_0)\cdot\frac{\e{\te}\cdot\hat\mu(\om_0-1)}{f(\om_{-1}+1)\hat\mu(\om_0)}-p(\om_{-1},\,\om_0+1)\\
% &\quad-p(\om_{-1},\,\om_0)+p(\om_{-1},\,\om_0+1)\\
% &\quad+p(\om_0+1,\,\om_1-1)\cdot\frac{\hat\mu(\om_0+1)\cdot\e{-\si}\cdot f(\om_1)}{\hat\mu(\om_0)}-p(\om_0,\,\om_1)\\
% &\quad-p(\om_0+1,\,\om_1)+p(\om_0,\,\om_1)\\
% &\quad+q(\om_{-1},\,\om_{a-1})-q(\om_{a-1},\,\om_{-1})\\
% &\quad+q(\om_{b+1},\,\om_1)-q(\om_1,\,\om_{b+1})\\
% &\quad+q(\om_{-1}-1,\,\om_0+1)\cdot\frac{\e{-\te}\cdot f(\om_{-1})\hat\mu(\om_0+1)}{\hat\mu(\om_0)}-q(\om_{-1},\,\om_0)\\
% &\quad-q(\om_{-1},\,\om_0+1)+q(\om_{-1},\,\om_0)\\
% &\quad+q(\om_0,\,\om_1+1)\cdot\frac{\hat\mu(\om_0-1)\cdot\e{\si}}{\hat\mu(\om_0)f(\om_1+1)}-q(\om_0+1,\,\om_1)\\
% &\quad-q(\om_0,\,\om_1)+q(\om_0+1,\,\om_1).
% \ea
% \]
% We continue applying \eqref{eq:telesc} on indices \(-1\), \(0\) and \(1\):
\be
\ba
A&=p(\om_{b+1},\,\om_{a-1})-p(\om_{a-1},\,\om_{b+1})+q(\om_{b+1},\,\om_{a-1})-q(\om_{a-1},\,\om_{b+1})\\
&\quad+p(\om_{-1}+1,\,\om_0)\cdot\frac{\e{\te}\cdot\hat\mu(\om_0-1)}{f(\om_{-1}+1)\hat\mu(\om_0)}-p(\om_0,\,\om_{-1})\\
&\quad+p(\om_0+1,\,\om_1-1)\cdot\frac{\hat\mu(\om_0+1)\cdot\e{-\si}\cdot f(\om_1)}{\hat\mu(\om_0)}-p(\om_0+1,\,\om_1)\\
&\quad+p(\om_0,\,\om_1)-p(\om_1,\,\om_0)\\
&\quad+q(\om_{-1}-1,\,\om_0+1)\cdot\frac{\e{-\te}\cdot f(\om_{-1})\hat\mu(\om_0+1)}{\hat\mu(\om_0)}-q(\om_{-1},\,\om_0+1)\\
&\quad+q(\om_{-1},\,\om_0)-q(\om_0,\,\om_{-1})\\
&\quad+q(\om_0,\,\om_1+1)\cdot\frac{\hat\mu(\om_0-1)\cdot\e{\si}}{\hat\mu(\om_0)f(\om_1+1)}-q(\om_1,\,\om_0).
\ea\label{eq:a}
\ee

Checking if random walking shocks with second class particles emerge in a model now simplifies to checking the existence of a measure \(\hat\mu\) with which the terms \(A\), \(B+D\), \(C+E\) each do not depend on \(\om_{-1}\), \(\om_0\), \(\om_1\). (Recall that \(\vp\) depends on configurations of sites \(a\dots b\), hence \(A\) is allowed to depend on \(\om_{a-1}\) and \(\om_{b+1}\).) When \(\omax\) is finite, one also has to take into account the fact that \(B+D\), \(C+E\) and \(A\) are taken under \(\Ev_{-1}\), \(\Ev_1\) and \(\Ev_0\), respectively, which give zero weight on \(\om_{-1}=\omax\), \(\om_1=\omax\) and \(\om_0=\omax\), respectively.

Substituting the rates of ASEP, \eqref{eq:aseprel} and \eqref{eq:asepmu} give for \eqref{eq:bd}, \eqref{eq:ce}:
\[
B+D=q\cdot\frac{1-\vr}{1-\la}=Q,\qquad C+E=p\cdot\frac{1-\la}{1-\vr}=P,
\]
in agreement with \eqref{eq:aseppq}. We obtain, also using \eqref{eq:asepro}, from \eqref{eq:a}
\[
A=(p-q)(\om_{b+1}-\om_{a-1})-1.
\]
Its expectation w.r.t.\ \(\Ev_0\) can be directly computed since \(\vp\) does not depend on \(\om_{b+1}\) nor \(\om_{a-1}\), and \(\Ev_0\) is product:
\[
\Ev_0A=(p-q)\cdot(\la-\vr)-1=-q\cdot\frac{1-\vr}{1-\la}-p\cdot\frac{1-\la}{1-\vr}
\]
which finishes the proof for ASEP by \eqref{eq:ae}, \eqref{eq:main} and \eqref{eq:aseppq}.

The rates \(p(y,\,z)=f(y)=\e{\beta(y-1/2)}\) and \(q(y,\,z)=0\) of the exponential GZRP result in
\[
B+D=0=Q,\qquad C+E=\e{\te}-\e{\si}=P,
\]
where we used \eqref{eq:om}, \eqref{eq:expmu}, \eqref{eq:expz} together with \eqref{eq:exprel}. This again agrees with \eqref{eq:blppq}. From \eqref{eq:a} we get
\[
A=f(\om_{b+1})-f(\om_{a-1}),
\]
which becomes \(\e{\si}-\e{\te}\) under the expectation of the product measure \(\un\nu_0\), also in agreement with \eqref{eq:blppq}. Similar computations for the exponential BLP lead to
\[
\begin{gathered}
B+D=\e{-\si}-\e{-\te}=Q,\qquad C+E=\e{\te}-\e{\si}=P,\\
A=f(\om_{b+1})+f(-\om_{a-1})-f(\om_{a-1})-f(-\om_{b+1}),
\end{gathered}
\]
which becomes \(\e{\si}+\e{-\te}-\e{\te}-\e{-\si}\) under \(\Ev_0\).
\end{proof}

\subsection{The proof for multiple shocks}

\begin{proof}[Proof of Theorem \ref{tm:mush}]
Expectation w.r.t.\ \(\un\nu^{\un\si,\un m}\) \eqref{eq:unnudef} will be denoted by \(\Ev^{\un\si,\un m}\). Take a finite cylinder function that depends on pairs \((\om_i,\,\ze_i)\) for \(a\le i\le b\) and assume that \(m_i=0\) for \(i\le a\) and \(i\ge b\). (This implies, in particular, a finite number of second class particles in the system.) We will take the \(\Ev^{\un\si,\un m}\) expectation of \eqref{eq:altcougen} and change variables in some terms so that after the change we have \(\vp(\un\om,\,\un\ze)\) in each one. This expectation involves summations for all \(\om_i\) variables. We shall change the summation variables in each term where the argument of \(\vp\) is modified such that, after the change, each term will contain \(\vp(\un\om,\,\un\ze)\). Lines \eqref{eq:altomb} - \eqref{eq:altomj} describe jumps of second class particles. For all the other lines, the second class particle stays at site \(i\), and the change of variables can be done without much complications. In the term with \(\vp(\un\om^{i,i+1},\,\un\ze^{i,i+1})\) in the first line, for example, we will put \(\un\om^{i,i+1}\) and \(\un\ze^{i,i+1}\) as our new summation variables. This step will change the arguments of the rates, and it will also bring in factors of \(\nu^{\si_i,m_i}\) and \(\nu^{\si_{i+1},m_{i+1}}\). We will also make use of \(\ze_i=\om_i+m_i\) which is a probability one event under \(\Ev^{\un\si,\un m}\).

For lines \eqref{eq:altomb} - \eqref{eq:altomj}, however, notice that the change of variables influences the number of second class particles, and the measures \(\un\nu^{\un\si,\un m}\) are singular to each other with different vectors \(\un m\). Therefore with the change of variables we will also include a replacement of measures accordingly. 

These changes result in 
\begin{multline}
\Ev^{\un\si,\un m}(L\vp)(\un\om,\,\un\ze)=\sum_{i=a-1}^b\biggl\{\Ev^{\un\si,\un m}\bigl\{\vp(\un\om,\,\un\ze)\cdot A_i\bigr\}\\
+\Ev^{\un\si^{i+1,i},\un m^{i+1,i}}\bigl\{\vp(\un\om,\,\un\ze)\cdot[B_i+D_i]\bigr\}+\Ev^{\un\si^{i,i+1},\un m^{i,i+1}}\bigl\{\vp(\un\om,\,\un\ze)\cdot[C_i+E_i]\bigr\}\biggr\},\label{eq:mul}
\end{multline}
where \(B_i\) comes from the term with \(\vp(\un\om^{i,i+1},\,\un\ze)\) of line \eqref{eq:altomb}, \(C_i\) comes from the term with \(\vp(\un\om,\,\un\ze^{i,i+1})\) of line \eqref{eq:altzej}, \(D_i\) comes from the term with \(\vp(\un\om,\,\un\ze^{i+1,i})\) of line \eqref{eq:altzeb}, \(E_i\) comes from the term with \(\vp(\un\om^{i+1,i},\,\un\ze)\) of line \eqref{eq:altomj}, and \(A_i\) from all other terms. The aim is again to show that \(B_i+D_i\), \(C_i+E_i\), and the sum of \(A_i\)'s each do not depend on \(\om_a\dots \om_b\). In fact, they will be identified as the rates \eqref{eq:asepshrj}, \eqref{eq:asepshrb}, and \eqref{eq:zrpshrj} -- \eqref{eq:blpshrb}, respectively. We now compute the quantities \(A_i\) - \(E_i\).
\be
\ba
A_i&=p(\om_i+1,\,\om_{i+1}+m_{i+1}-1)\cdot\frac{\hat\mu^{\si_i,m_i}(\om_i+1)\hat\mu^{\si_{i+1},m_{i+1}}(\om_{i+1}-1)}{\hat\mu^{\si_i,m_i}(\om_i)\hat\mu^{\si_{i+1},m_{i+1}}(\om_{i+1})}\\
&\quad-p(\om_i,\,\om_{i+1}+m_{i+1})\\
&\quad-p(\om_i,\,\om_{i+1})+p(\om_i,\,\om_{i+1}+m_{i+1})\\
&\quad-p(\om_i+m_i,\,\om_{i+1}+m_{i+1})+p(\om_i,\,\om_{i+1}+m_{i+1})\\
&\quad+q(\om_i+m_i-1,\,\om_{i+1}+1)\cdot\frac{\hat\mu^{\si_i,m_i}(\om_i-1)\hat\mu^{\si_{i+1},m_{i+1}}(\om_{i+1}+1)}{\hat\mu^{\si_i,m_i}(\om_i)\hat\mu^{\si_{i+1},m_{i+1}}(\om_{i+1})}\\
&\quad-q(\om_i+m_i,\,\om_{i+1})\\
&\quad-q(\om_i+m_i,\,\om_{i+1}+m_{i+1})+q(\om_i+m_i,\,\om_{i+1})\\
&\quad-q(\om_i,\,\om_{i+1})+q(\om_i+m_i,\,\om_{i+1})\\
&=p(\om_i+1,\,\om_{i+1}+m_{i+1}-1)\cdot\frac{\hat\mu^{\si_i,m_i}(\om_i+1)\hat\mu^{\si_{i+1},m_{i+1}}(\om_{i+1}-1)}{\hat\mu^{\si_i,m_i}(\om_i)\hat\mu^{\si_{i+1},m_{i+1}}(\om_{i+1})}\\
&\quad-p(\om_i,\,\om_{i+1})+p(\om_i,\,\om_{i+1}+m_{i+1})-p(\om_i+m_i,\,\om_{i+1}+m_{i+1})\\
&\quad+q(\om_i+m_i-1,\,\om_{i+1}+1)\cdot\frac{\hat\mu^{\si_i,m_i}(\om_i-1)\hat\mu^{\si_{i+1},m_{i+1}}(\om_{i+1}+1)}{\hat\mu^{\si_i,m_i}(\om_i)\hat\mu^{\si_{i+1},m_{i+1}}(\om_{i+1})}\\
&\quad-q(\om_i+m_i,\,\om_{i+1}+m_{i+1})+q(\om_i+m_i,\,\om_{i+1})-q(\om_i,\,\om_{i+1}),
\ea\label{eq:ai}
\ee
\be
\ba
B_i+D_i&=\bigl[p(\om_i+1,\,\om_{i+1}-1)-p(\om_i+1,\,\om_{i+1}+m_{i+1}-1)\bigr]\\
&\qquad\times\frac{\hat\mu^{\si_i,m_i}(\om_i+1)\hat\mu^{\si_{i+1},m_{i+1}}(\om_{i+1}-1)}{\hat\mu^{\si^{i+1,i}_i,m_i+1}(\om_i)\hat\mu^{\si^{i+1,i}_{i+1},m_{i+1}-1}(\om_{i+1})}\\
&+\bigl[q(\om_i+m_i,\,\om_{i+1}+m_{i+1})-q(\om_i+m_i,\,\om_{i+1})\bigr]\\
&\qquad\times\frac{\hat\mu^{\si_i,m_i}(\om_i)\hat\mu^{\si_{i+1},m_{i+1}}(\om_{i+1})}{\hat\mu^{\si^{i+1,i}_i,m_i+1}(\om_i)\hat\mu^{\si^{i+1,i}_{i+1},m_{i+1}-1}(\om_{i+1})},
\ea\label{eq:bdi}
\ee
\be
\ba
C_i+E_i&=\bigl[p(\om_i+m_i,\,\om_{i+1}+m_{i+1})-p(\om_i,\,\om_{i+1}+m_{i+1})\bigr]\\
&\qquad\times\frac{\hat\mu^{\si_i,m_i}(\om_i)\hat\mu^{\si_{i+1},m_{i+1}}(\om_{i+1})}{\hat\mu^{\si^{i,i+1}_i,m_i-1}(\om_i)\hat\mu^{\si^{i,i+1}_{i+1},m_{i+1}+1}(\om_{i+1})}\\
&+\bigl[q(\om_i-1,\,\om_{i+1}+1)-q(\om_i+m_i-1,\,\om_{i+1}+1)\bigr]\\
&\qquad\times\frac{\hat\mu^{\si_i,m_i}(\om_i-1)\hat\mu^{\si_{i+1},m_{i+1}}(\om_{i+1}+1)}{\hat\mu^{\si^{i,i+1}_i,m_i-1}(\om_i)\hat\mu^{\si^{i,i+1}_{i+1},m_{i+1}+1}(\om_{i+1})}.
\ea\label{eq:cei}
\ee
While these formulas are much more general and complicated than \eqref{eq:a}, \eqref{eq:bd} and \eqref{eq:ce}, they can be used to verify Theorem \ref{tm:mush}.

We start with plugging in the rates of ASEP and \eqref{eq:asepsokfelt}, \eqref{eq:asepsokmuhat}.
\[
\ba
A_i&=p{\bf1}\{\om_i=0,\,\om_{i+1}=1,\,m_i=0,\,m_{i+1}=0\}\cdot\frac{\vr_i(1-\vr_{i+1})}{(1-\vr_i)\vr_{i+1}}\\
&\quad-p{\bf1}\{\om_i=1,\,\om_{i+1}=0\}+p{\bf1}\{\om_i=1,\,\om_{i+1}=0,\,m_{i+1}=0\}\\
&\quad-p{\bf1}\{\om_i=0,\,\om_{i+1}=0,\,m_i=1,\,m_{i+1}=0\}\\
&\quad-p{\bf1}\{\om_i=1,\,\om_{i+1}=0,\,m_i=0,\,m_{i+1}=0\}\\
&\quad+q{\bf1}\{\om_i=1,\,\om_{i+1}=0,\,m_i=0,\,m_{i+1}=0\}\cdot\frac{(1-\vr_i)\vr_{i+1}}{\vr_i(1-\vr_{i+1})}\\
&\quad-q{\bf1}\{\om_i=0,\,\om_{i+1}=1,\,m_i=0,\,m_{i+1}=0\}\\
&\quad-q{\bf1}\{\om_i=0,\,\om_{i+1}=0,\,m_i=0,\,m_{i+1}=1\}\\
&\quad+q{\bf1}\{\om_i=0,\,\om_{i+1}=1,\,m_i=0\}-q{\bf1}\{\om_i=0,\,\om_{i+1}=1\}\\
%&\ba
%={\bf1}\{m_i=0,\,m_{i+1}=0\}\Bigl\{&{\bf1}\{\om_i=0,\,\om_{i+1}=1\}\cdot\Bigl[p\frac{\vr_i(1-\vr_{i+1})}{(1-\vr_i)\vr_{i+1}}-q\Bigr]\\
%+&{\bf1}\{\om_i=1,\,\om_{i+1}=0\}\cdot\Bigl[q\frac{(1-\vr_i)\vr_{i+1}}{\vr_i(1-\vr_{i+1})}-p\Bigr]\Bigr\}
%\ea\\
%&-{\bf1}\{m_i=0,\,m_{i+1}=1\}\bigl[p{\bf1}\{\om_i=1\}+q{\bf1}\{\om_i=0\}\bigr]\\
%&-{\bf1}\{m_i=1,\,m_{i+1}=0\}\bigl[p{\bf1}\{\om_{i+1}=0\}+q{\bf1}\{\om_{i+1}=1\}\bigr]\\
%&={\bf1}\{m_i=0,\,m_{i+1}=0\}\bigl[{\bf1}\{\om_i=0,\,\om_{i+1}=1\}-{\bf1}\{\om_i=1,\,\om_{i+1}=0\}\bigr](p-q)\\
%&-{\bf1}\{m_i=0,\,m_{i+1}=1\}\bigl[p{\bf1}\{\om_i=1\}+q{\bf1}\{\om_i=0\}\bigr]\\
%&-{\bf1}\{m_i=1,\,m_{i+1}=0\}\bigl[p{\bf1}\{\om_{i+1}=0\}+q{\bf1}\{\om_{i+1}=1\}\bigr]\\
&=(p-q)\bigl[\om_{i+1}(1-m_{i+1})-\om_i(1-m_i)\bigr]\\
&\quad-pm_i(1-m_{i+1})-q(1-m_i)m_{i+1}.
%%&=p{\bf1}\{\om_i=0,\,\om_{i+1}=1,\,m_i=0,\,m_{i+1}=0\}\cdot\frac{\vr_i(1-\vr_{i+1})}{(1-\vr_i)\vr_{i+1}}\\
%%&-p{\bf1}\{\om_i=1,\,\om_{i+1}=0,\,m_i=0\}\\
%%&-p{\bf1}\{\om_i=0,\,\om_{i+1}=0,\,m_i=1,\,m_{i+1}=0\}\\
%%&+q{\bf1}\{\om_i=1,\,\om_{i+1}=0,\,m_i=0,\,m_{i+1}=0\}\cdot\frac{(1-\vr_i)\vr_{i+1}}{\vr_i(1-\vr_{i+1})}\\
%%&-q{\bf1}\{\om_i=0,\,\om_{i+1}=0,\,m_i=0,\,m_{i+1}=1\}\\
%%&-q{\bf1}\{\om_i=0,\,\om_{i+1}=1,\,m_{i+1}=0\}
\ea
\]
In the second equality we used that \(\hat\mu^{\vr,1}(1)=0\) in \(\Ev^{\un\si,\un m}\), hence any situation with \(\om_j=1\), \(m_j=1\) does not occur. Validity of the equation can be directly checked for the remaining 9 cases of \(\om's\) and \(m\)'s being 0 or 1.
\[
\ba
B_i+D_i&=p{\bf1}\{\om_i=0,\,\om_{i+1}=1,\,m_i=0,\,m_{i+1}=1\}\frac{\vr_i}{\vr_{i+1}}\\
&+q{\bf1}\{\om_i=0,\,\om_{i+1}=0,\,m_i=0,\,m_{i+1}=1\}\frac{1-\vr_i}{1-\vr_{i+1}}\\
&=p(1-m_i)m_{i+1}\cdot\frac{\vr_i}{\vr_{i+1}}=(1-m_i)m_{i+1}\cdot\bigl[(1-\vr_i)q+\vr_ip\bigr]
\ea
\]
using \eqref{eq:asepsokfelt} and the fact that the indicators imply \(m^{i+1,i}_i=1\) and therefore necessarily \(\om_i=0\) under \(\Ev^{\un\si^{i+1,i},\un m^{i+1,i}}\). Similarly,
\[
\ba
C_i+E_i&=p{\bf1}\{\om_i=0,\,\om_{i+1}=0,\,m_i=1,\,m_{i+1}=0\}\frac{1-\vr_{i+1}}{1-\vr_{i-1}}\\
&+q{\bf1}\{\om_i=1,\,\om_{i+1}=0,\,m_i=1,\,m_{i+1}=0\}\frac{\vr_{i+1}}{\vr_{i-1}}\\
&=qm_i(1-m_{i+1})\cdot\frac{\vr_{i+1}}{\vr_{i-1}}=m_i(1-m_{i+1})\cdot\bigl[(1-\vr_{i+1})p+\vr_{i+1}q\bigr].
\ea
\]
Compare these with \eqref{eq:asepshrb} and \eqref{eq:asepshrj}. Now we redistribute some terms from the \(A_i\)'s and write \eqref{eq:mul} as
\begin{multline*}
\Ev^{\un\si,\un m}(L\vp)(\un\om,\,\un\ze)\\
\ba
=\sum_{i=a-1}^b\Bigl\{&(1-m_i)m_{i+1}\bigl[(1-\vr_i)q+\vr_ip\bigr]\\
&\qquad\times\bigl[\Ev^{\un\si^{i+1,i},\un m^{i+1,i}}\vp(\un\om,\,\un\ze)-\Ev^{\un\si,\un m}\vp(\un\om,\,\un\ze)\bigr]\\
+\,&m_i(1-m_{i+1})\bigl[(1-\vr_{i+1})p+\vr_{i+1}q\bigr]\\
&\qquad\times\bigl[\Ev^{\un\si^{i,i+1},\un m^{i,i+1}}\vp(\un\om,\,\un\ze)-\Ev^{\un\si,\un m}\vp(\un\om,\,\un\ze)\bigr]\Bigr\}\\
+\,&(p-q)\Ev^{\un\si,\un m}\vp(\un\om,\,\un\ze)\sum_{i=a-1}^b\bigl[\om_{i+1}(1-m_{i+1})-\om_i(1-m_i)\\
&\qquad-\vr_{i+1}m_i(1-m_{i+1})+\vr_i(1-m_i)m_{i+1}\bigr].
\ea
\end{multline*}
The proof for ASEP is done as soon as we see that the last two lines sum up to zero. Without loss of generality, let \(m_i=0\) for each \(i\le a\) and \(i\ge b\). The terms with \(\om\)'s are telescopic, and the terms in the bracket will only have \(\om_{a-1}\) and \(\om_{b+1}\) in them. These become \(\vr_{a-1}\) and \(\vr_{b+1}\), respectively, under the product expectation. Therefore the expectation transforms the summation into
\be
\vr_{b+1}-\vr_{a-1}-\sum_{i=a-1}^b\bigl[\vr_{i+1}m_i(1-m_{i+1})-\vr_i(1-m_i)m_{i+1}\bigr].\label{eq:asepais}
\ee
Since there are finitely many second class particles, there is a leftmost site \(r_1\) such that \(m_{r_1}=0\) and \(m_{r_1+1}=1\). We successively define
\[
\ba
\ell_k:&=\inf\{i\,:\,i>r_k,\ m_i=1,\ m_{i+1}=0\}\text{, and}\\
r_{k+1}:&=\inf\{i\,:\,i>r_k,\ m_i=0,\ m_{i+1}=1\}.
\ea
\]
We let \(\inf\emptyset=\infty\), and \(K=\max\{k\,:\,r_k<\infty\}\). For any \(i\) with \(\ell_k<i<r_{k+1}\) we have \(m_i=m_{i+1}=0\), and for any \(i\) with \(r_k<i<\ell_k\) we have \(m_i=m_{i+1}=1\). Therefore, the summand is zero for all these cases, and
\[
\ba
&\sum_{i=a-1}^b\bigl[\vr_{i+1}m_i(1-m_{i+1})-\vr_i(1-m_i)m_{i+1}\bigr]\\
&=\sum_{k=1}^K[\vr_{\ell_k+1}-\vr_{r_k}]=\vr_{\ell_K+1}-\vr_{r_1}+\sum_{k=1}^{K-1}[\vr_{\ell_k+1}-\vr_{r_{k+1}}].
\ea
\]
Notice also that by \eqref{eq:asepsokfelt} \(\vr_i\) is unchanged over an interval without second class particles. Hence the sum in the last display is zero, and \(\vr_{\ell_K+1}=\vr_{b+1}\), \(\vr_{r_1}=\vr_{a-1}\) implies that \eqref{eq:asepais} is zero which finishes the proof for ASEP.

We now turn to the exponential GZRP and BLP models. Substituting into \eqref{eq:ai}, \eqref{eq:bdi} and \eqref{eq:cei} gives
\[
\ba
A_i&=\e{\beta(\om_{i+1}+m_{i+1})-\beta/2}-\e{\beta(\om_i+m_i)-\beta/2},\\
B_i+D_i&=0,\\
C_i+E_i&=\e{\si_i+\beta m_i}-\e{\si_i}
\ea
\]
for the GZRP, and
\[
\ba
A_i&=\bigl(\e{\beta(\om_{i+1}+m_{i+1})-\beta/2}-\e{\beta(\om_i+m_i)-\beta/2}\bigr)-\bigl(\e{-\beta\om_{i+1}-\beta/2}-\e{-\beta\om_i-\beta/2}\bigr),\\
B_i+D_i&=\e{-\si_{i+1}}-\e{-\si_{i+1}-\beta m_{i+1}},\\
C_i+E_i&=\e{\si_i+\beta m_i}-\e{\si_i}
\ea
\]
for the BLP, where we also used \eqref{eq:expz}. Thus we identified these terms with \eqref{eq:zrpshrj} -- \eqref{eq:blpshrb} as required. Notice that \(A_i\) is of gradient form, thus summing it in the first line of \eqref{eq:mul} gives
\begin{multline*}
\sum_{i=a-1}^bA_i\\
\ba
&=\e{\beta(\om_{b+1}+m_{b+1})-\beta/2}-\e{\beta(\om_{a-1}+m_{a-1})-\beta/2}\\
&=f(\om_{b+1}+m_{b+1})-f(\om_{a-1}+m_{a-1})\qquad\text{for GZRP,}\\
&=\bigl(\e{\beta(\om_{b+1}+m_{b+1})-\beta/2}-\e{\beta(\om_{a-1}+m_{a-1})-\beta/2}\bigr)-\bigl(\e{-\beta\om_{b+1}-\beta/2}-\e{-\beta\om_{a-1}-\beta/2}\bigr)\\
&=f(\om_{b+1}+m_{b+1})-f(\om_{a-1}+m_{a-1})-f(-\om_{b+1})+f(-\om_{a-1})\qquad\text{for BLP.}
\ea
\end{multline*}
We again can assume finitely many second class particles inside the discrete interval \([a,\,b]\). Product structure of the measure then allows us to compute the above and show that the first line of \eqref{eq:mul} equals
\begin{multline*}
\sum_{i=a-1}^b\Ev^{\un\si,\un m}\bigl\{\vp(\un\om,\,\un\ze)\cdot A_i\bigr\}\\
\ba
&=\bigl[\e{\si_{b+1}}-\e{\si_{a-1}}\bigr]\cdot\Ev^{\un\si,\un m}\vp(\un\om,\,\un\ze)&&\quad\text{for GZRP,}\\
&=\bigl[\e{\si_{b+1}}-\e{\si_{a-1}}-\e{-\si_{b+1}}+\e{-\si_{a-1}}\bigr]\cdot\Ev^{\un\si,\un m}\vp(\un\om,\,\un\ze)&&\quad\text{for BLP.}
\ea
\end{multline*}
These brackets can be broken up into (notice that \(\si_i\) is constant near the boundaries)
\[
\sum_{i=a-1}^b\bigl[\e{\si_{i+1}}-\e{\si_i}\bigr]=\sum_{i=a-1}^b\bigl[\e{\si_i}-\e{\si_{i-1}}\bigr]=\sum_{i=a-1}^b\bigl[\e{\si_i}-\e{\si_i+\beta m_i}\bigr]
\]
for the GZRP, and
\[
\sum_{i=a-1}^b\bigl[\e{\si_{i+1}}-\e{\si_i}-\e{-\si_{i+1}}+\e{-\si_i}\bigr]=\sum_{i=a-1}^b\bigl[\e{\si_i}-\e{\si_i+\beta m_i}-\e{-\si_{i+1}}+\e{-\si_{i+1}-\beta m_{i+1}}\bigr]
\]
for the BLP. These give the negative terms in \eqref{eq:mush} and the proof is done.
\end{proof}

\section{Branching coalescing random walk}\label{sc:bcrw}

\subsection{The model}

In this section we investigate the random walking shock structure found in the branching coalescing random walk model (BCRW) by Krebs, Jafarpour, Sch\"utz \cite{bcrw}. This model does not conserve the particle number, it is not a member of the family of Section \ref{sc:fam}. However, the result we show below is very similar to Theorem \ref{tm:main}.

Besides \eqref{eq:ugrdef}, we introduce \(\un\om^{i\up}\) and \(\un\om^{i\dn}\) corresponding to particle creation and annihilation on site $i$ by
\[
\om^{i\up}_k=\left\{\ba
&\om_k,&&\text{for }k\ne i,\\
&\om_k+1,&&\text{for }k=i,
\ea\right.\qquad
\om^{i\dn}_k=\left\{\ba
&\om_k,&&\text{for }k\ne i,\\
&\om_k-1,&&\text{for }k=i.
\ea\right.
\]
The BCRW is a particle system on \(\{0,\,1\}^\Zb\), and the dynamics consists of three types of processes. Asymmetric nearest neighbour jumps (like in the ASEP):
\[
\ba
\un\om&\to\un\om^{i,i+1}&\qquad&\text{with rate }p\cdot\om_i(1-\om_{i+1}),\\
\un\om&\to\un\om^{i+1,i}&\qquad&\text{with rate }q\cdot(1-\om_i)\om_{i+1},
\ea
\]
coalescence (i.e., merging of two particles into one) from the left and from the right:
\[
\ba
\un\om&\to\un\om^{i\dn}&\qquad&\text{with rate }c_r\cdot\om_i\om_{i+1},\\
\un\om&\to\un\om^{i+1\dn}&\qquad&\text{with rate }c_l\cdot\om_i\om_{i+1},
\ea
\]
and branching (i.e., creation of a new particle next to an existing one) to the left and to the right:
\[
\ba
\un\om&\to\un\om^{i\up}&\qquad&\text{with rate }b_l\cdot(1-\om_i)\om_{i+1},\\
\un\om&\to\un\om^{i+1\up}&\qquad&\text{with rate }b_r\cdot\om_i(1-\om_{i+1}).
\ea
\]
The factors in the rates here are positive real numbers. The generator of the process (with a finite cylinder function \(\vp\)) is
\begin{multline}
(L\vp)(\un\om)=\\
\ba
\sum_{i=a-1}^b\bigl\{&p\cdot\om_i(1-\om_{i+1})[\vp(\un\om^{i,i+1})-\vp(\un\om)]+q\cdot(1-\om_i)\om_{i+1}[\vp(\un\om^{i+1,i})-\vp(\un\om)]\\
+\,&c_r\cdot\om_i\om_{i+1}[\vp(\un\om^{i\dn})-\vp(\un\om)]+c_l\cdot\om_i\om_{i+1}[\vp(\un\om^{i+1\dn})-\vp(\un\om)]\\
+\,&b_l\cdot(1-\om_i)\om_{i+1}[\vp(\un\om^{i\up})-\vp(\un\om)]+b_r\cdot\om_i(1-\om_{i+1})[\vp(\un\om^{i+1\up})-\vp(\un\om)]\bigl\}.
\ea\label{eq:bcrwg}
\end{multline}

We will abbreviate \(B=b_l+b_r\) and \(C=c_l+c_r\). Define, for \(\vr\in[0,\,1]\), the one-site marginal \(\mu^\vr\) by
\[
\mu^\vr(1)=1-\mu^\vr(0)=\vr.
\]
It is known (and will be apparent from our computations as well) that the Bernoulli product distribution
\be
\bigotimes_{i\in\Zb}\mu^{\vr^*}\text{ with the specific density value } \vr^*=\frac B{B+C}\label{eq:rocs}
\ee
is a translation-invariant stationary distribution. Besides this, there is a trivial stationary measure, which is the totally empty lattice with probability one: $\bigotimes_{i\in\Zb}\mu^0$.

\subsection{Shock measure}\label{sc:bcrwres}

In \cite{bcrw} the following shock measure was considered (adapted to infinite volume):
\be
\bigotimes_{i\le j}\mu^{\vr^*}\bigotimes_{i>j}\mu^0.\label{eq:gumes}
\ee
This measure is the \(\vr^*\)-stationary distribution on the left of \(j\), and forces the empty configuration on the right of \(j\). It has been shown that this shock structure performs a biased simple random walk when the algebraic relations \eqref{eq:bcrwfelt} below hold for the rates.

This result is similar to that for the ASEP \cite{qse} and for the BLP \cite{sokvalak}. An important difference is however, that whereas in the ASEP and the BLP a good definition for the microscopic position of the shock is the position of the second class particle, in the BCRW this is not possible (the basic coupling does not conserve the number of second class particles). Here, an obvious definition of the shock position is the position of the rightmost particle. Indeed, using the result of \cite{bcrw}, it can easily be shown that the stationary distribution seen from the rightmost particle is a product measure with density $\vr^*$ on negative sites and density zero on positive ones, which is a similar result to that of \cite{dls} and \cite{valak}. This naturally raises the question whether the rightmost particle itself performs a simple random walk in the BCRW. 

Following the previous strategy we define the shock measure
\[
\un\mu_j=\bigotimes_{i<j}\mu^{\vr^*}\bigotimes_{i=j}\mu^1\bigotimes_{i>j}\mu^0,
\]
which is similar to the previous shock measure except that at site \(j\) we force the presence of a particle. This one is of course the rightmost particle in the system.
\begin{tm}\label{tm:bcrw}
The identity
\be
\frac{\di}{\di t}\un\mu_jS(t)\Bigr|_{t=0}=P\cdot[\un\mu_{j+1}-\un\mu_j]+Q\cdot[\un\mu_{j-1}-\un\mu_j]\label{eq:bcrw}
\ee
holds if
\be
p=b_r\cdot\frac CB.\label{eq:bcrwfelt}
\ee
In this case
\be
P=p\cdot\frac{C+B}C,\qquad\text{and}\qquad Q=q\cdot\frac C{C+B}+c_l\cdot\frac B{C+B},\label{eq:bcrwpq}
\ee
which shows that the rightmost particle performs a simple random walk.
\end{tm}
Notice that \eqref{eq:bcrwfelt} does not depend on \(q\), the jump rate of particles away from the zero-density region.
We remark that both the result in \cite{bcrw} and this theorem has a mirror-sym\-met\-ric version with the empty configuration on the left and density \(\vr^*\) on the right, \(p\), \(c_r\), and \(b_r\) interchanged with \(q\), \(c_l\), and \(b_l\), respectively.

Notice also that
\[
\sum_{k=-\infty}^j(1-\vr^*)^{j-k}\vr^*\cdot\un\mu_k
\]
exactly equals \eqref{eq:gumes} (just condition on where the rightmost particle is according to \eqref{eq:gumes}). Due to linearity of both sides of \eqref{eq:bcrw}, our result implies the same equation for the measure \eqref{eq:gumes}, thus the result of \cite{bcrw} is recovered.

\subsection{Proof for BCRW}

\begin{proof}[Proof of Theorem \ref{tm:bcrw}]
Again we set \(j=0\). The idea is now familiar: first take a cylinder function of values \(\om_a\dots\om_b\) with \(a<-1<1<b\) and apply the expectation \(\Ev_0\) w.r.t.\ \(\un\mu_0\) of \eqref{eq:bcrwg}. Making use of \(\om_0=1\) and \(\om_i=0\) for \(i>0\), this leads to
\[
\ba
\Ev_0(L\vp)(\un\om)=\Ev_0\Bigl\{&\sum_{i=a-1}^{-2}\bigl\{p\cdot\om_i(1-\om_{i+1})[\vp(\un\om^{i,i+1})-\vp(\un\om)]\\
&\qquad+q\cdot(1-\om_i)\om_{i+1}[\vp(\un\om^{i+1,i})-\vp(\un\om)]\\
&\qquad+c_r\cdot\om_i\om_{i+1}[\vp(\un\om^{i\dn})-\vp(\un\om)]\\
&\qquad+c_l\cdot\om_i\om_{i+1}[\vp(\un\om^{i+1\dn})-\vp(\un\om)]\\
&\qquad+b_l\cdot(1-\om_i)\om_{i+1}[\vp(\un\om^{i\up})-\vp(\un\om)]\\
&\qquad+b_r\cdot\om_i(1-\om_{i+1})[\vp(\un\om^{i+1\up})-\vp(\un\om)]\bigl\}\\
&+q\cdot(1-\om_{-1})[\vp(\un\om^{0,-1})-\vp(\un\om)]\\
&+c_r\cdot\om_{-1}[\vp(\un\om^{-1\dn})-\vp(\un\om)]\\
&+c_l\cdot\om_{-1}[\vp(\un\om^{0\dn})-\vp(\un\om)]\\
&+b_l\cdot(1-\om_{-1})[\vp(\un\om^{-1\up})-\vp(\un\om)]\\
&+p\cdot[\vp(\un\om^{0,1})-\vp(\un\om)]\\
&+b_r\cdot[\vp(\un\om^{1\up})-\vp(\un\om)]\Bigr\}.
\ea
\]
Then change variables to restore \(\vp(\un\om)\) in each term, and watch for singularity of \(\Ev_j\)'s for different configurations. In fact \(\un\om^{0,1}\) and \(\un\om^{1\up}\) will be taken under \(\Ev_1\), \(\un\om^{0,-1}\) and \(\un\om^{0\dn}\) will be taken under \(\Ev_{-1}\) after the change of variables. We arrive to
\be
\Ev_0(L\vp)(\un\om)=\Ev_0\Bigl\{\vp(\un\om)\Bigl[\sum_{i=a-1}^{-2}A_i+A\Bigr]\Bigr\}+\Ev_{-1}\bigl\{\vp(\un\om)\cdot D\bigr\}+\Ev_1\bigl\{\vp(\un\om)\cdot E\bigr\},\label{eq:bcrwe}
\ee
where
\[
\ba
A_i&=p\cdot(1-\om_i)\om_{i+1}-p\cdot\om_i(1-\om_{i+1})+q\cdot\om_i(1-\om_{i+1})-q\cdot(1-\om_i)\om_{i+1}\\
&+c_r\cdot(1-\om_i)\om_{i+1}\cdot\frac{\mu^{\vr^*}(1)}{\mu^{\vr^*}(0)}-c_r\cdot\om_i\om_{i+1}\\
&+c_l\cdot\om_i(1-\om_{i+1})\cdot\frac{\mu^{\vr^*}(1)}{\mu^{\vr^*}(0)}-c_l\cdot\om_i\om_{i+1}\\
&+b_l\cdot\om_i\om_{i+1}\cdot\frac{\mu^{\vr^*}(0)}{\mu^{\vr^*}(1)}-b_l\cdot(1-\om_i)\om_{i+1}\\
&+b_r\cdot\om_i\om_{i+1}\cdot\frac{\mu^{\vr^*}(0)}{\mu^{\vr^*}(1)}-b_r\cdot\om_i(1-\om_{i+1}),
\ea
\]
\[
\ba
A&=-q\cdot(1-\om_{-1})+c_r\cdot(1-\om_{-1})\cdot\frac{\mu^{\vr^*}(1)}{\mu^{\vr^*}(0)}-c_r\cdot\om_{-1}-c_l\cdot\om_{-1}\\
&+b_l\cdot\om_{-1}\frac{\mu^{\vr^*}(0)}{\mu^{\vr^*}(1)}-b_l\cdot(1-\om_{-1})-p-b_r,\\
D&=q\cdot\frac{\mu^{\vr^*}(0)}1+c_l\cdot\frac{\mu^{\vr^*}(1)}1,\\
E&=p\cdot(1-\om_0)\frac1{\mu^{\vr^*}(0)}+b_r\cdot\om_0\frac1{\mu^{\vr^*}(1)}.
\ea
\]
In the above we freely added or removed the a.s.\ one factors \(\om_j\) and \(1-\om_{j+1}\) for terms under \(\Ev_j\). The aim is again to show that \(D\) and \(E\) become independent of the variables and equal to \(Q\) and \(P\), respectively, while the \(A\) terms amount to \(-Q-P\) of \eqref{eq:bcrwpq}.
Plugging in the specific form \eqref{eq:rocs} of \(\vr^*\) simplifies things to
\[
\ba
A_i&=\bigl[p-q+c_r\cdot\frac BC-b_l\bigr]\cdot(\om_{i+1}-\om_i),\\
A&=\bigl[q-c_r\cdot\frac BC-C+b_l\cdot\frac CB+b_l\bigr]\cdot\om_{-1}-q+c_r\cdot\frac BC-b_l-p-b_r,\\
D&=q\cdot\frac C{B+C}+c_l\cdot\frac B{B+C},\\
E&=\bigl[-p-p\cdot\frac BC+b_r+b_r\cdot\frac CB\bigr]\cdot\om_0+p+p\cdot\frac BC.
\ea
\]
Since \(A_i\) is of gradient form, we can compute the sum in \eqref{eq:bcrwe}. We also use the fact that \(\vp\) does not depend on \(\om_{a-1}\), hence this latter has expectation \(\vr^*=B/[B+C]\). Then \eqref{eq:bcrwe} reads
\[
\ba
\Ev_0(L\vp)(\un\om)&=\Ev_0\vp(\un\om)\bigl\{\bigl[p-b_r\cdot\frac CB\bigr]\om_{-1}\\
&\qquad-p\cdot\bigl[\frac B{B+C}+1\bigr]-q\cdot\frac C{B+C}+\frac{b_rc_r-b_lc_l}{B+C}-b_r\bigr\}\\
&+\Ev_{-1}\vp(\un\om)\bigl\{q\cdot\frac C{B+C}+c_l\cdot\frac B{B+C}\bigr\}\\
&+\Ev_1\vp(\un\om)\bigl\{\bigl[-p-p\cdot\frac BC+b_r+b_r\cdot\frac CB\bigr]\cdot\om_0+p+p\cdot\frac BC\bigr\}.
\ea
\]
Substituting \eqref{eq:bcrwfelt} cancels the \(\om\)-dependence and results in
\[
\ba
\Ev_0(L\vp)(\un\om)&=\Ev_0\vp(\un\om)\bigl\{-q\cdot\frac C{B+C}-c_l\cdot\frac B{B+C}-p\cdot\frac{B+C}C\bigr\}\\
&+\Ev_{-1}\vp(\un\om)\bigl\{q\cdot\frac C{B+C}+c_l\cdot\frac B{B+C}\bigr\}\\
&+\Ev_1\vp(\un\om)\bigl\{p\cdot\frac{B+C}C\bigr\}.
\ea
\]
\end{proof}

\section*{Acknowledgments}

We thank Gunter Sch\"utz and Timo Sepp\"al\"ainen for fruitful discussions on the subject. We also thank our anonymous referees for valuable comments and suggestions. 

M.\ Bal\'azs was partially supported by the Hungarian Scientific Research Fund (OTKA) grants K-60708, F-67729, by the Bolyai Scholarship of the Hungarian Academy of Sciences, and by Morgan Stanley Mathematical Modeling Center.
A.\ R\'akos acknowledges financial support of the Hungarian Scientific Research Fund (OTKA) grants PD-72604, PD-78433 and from the Bolyai Scholarship of the Hungarian Academy of Sciences.

\appendix
\appendixpage

\section{Shock velocity for multiple shocks in the\\
ASEP}

Consider a multiple $n$-shock in the ASEP. We use a slightly modified notation: enumerate the second class particles \(k=1\dots n\) in spatial order, and denote by \(\vr_{(k)}\) the constant density value from the position of the \(k^\text{th}\) second class particle to the left neighbor of the site of the \(k+1^\text{st}\) second class particle (to infinity for \(k=n\)). We also denote by \(\vr_{(0)}\) the density value on the left of the \(1^\text{st}\) second class particle. This notation simplifies the arguments to come, and it is possible due to the specific fact that second class particles exclude each other and keep their order in ASEP. Condition \eqref{eq:asepsokfelt} then becomes
\[
\frac{\vr_{(k)}(1-\vr_{(k-1)})}{\vr_{(k-1)}(1-\vr_{(k)})}=\frac pq 
\]
in this notation, and jump rates \eqref{eq:asepshrj}, \eqref{eq:asepshrb} turn into
\begin{equation}\label{jumprates}
P_{(k)}=\frac{1-\vr_{(k)}}{1-\vr_{(k-1)}}\cdot p=\frac{\vr_{(k)}}{\vr_{(k-1)}}\cdot q\qquad\text{and}\qquad Q_{(k)}=\frac{1-\vr_{(k-1)}}{1-\vr_{(k)}}\cdot q=\frac{\vr_{(k-1)}}{\vr_{(k)}}\cdot p
\end{equation}
for the \(k^\text{th}\) second class particle. 

\begin{figure}
\begin{center}
\scalebox{0.99} % Change this value to rescale the drawing.
{
\begin{pspicture}(0,-2.7392187)(12.02,2.7792187)
\pscircle[linewidth=0.04,dimen=outer](1.5,1.5807812){0.3}
\psline[linewidth=0.04](0.0,1.2807813)(0.0,1.0807812)(12.0,1.0807812)(12.0,1.2807813)
\psline[linewidth=0.04cm](0.9984376,1.2807813)(0.9984376,1.0807812)
\psline[linewidth=0.04cm](1.9968752,1.2807813)(1.9968752,1.0807812)
\psline[linewidth=0.04cm](2.995313,1.2807813)(2.995313,1.0807812)
\psline[linewidth=0.04cm](3.9937503,1.2807813)(3.9937503,1.0807812)
\psline[linewidth=0.04cm](4.992188,1.2807813)(4.992188,1.0807812)
\psline[linewidth=0.04cm](5.990626,1.2807813)(5.990626,1.0807812)
\psline[linewidth=0.04cm](6.9890633,1.2807813)(6.9890633,1.0807812)
\psline[linewidth=0.04cm](7.9875007,1.2807813)(7.9875007,1.0807812)
\psline[linewidth=0.04cm](8.985938,1.2807813)(8.985938,1.0807812)
\psline[linewidth=0.04cm](9.984376,1.2807813)(9.984376,1.0807812)
\psline[linewidth=0.04cm](10.982814,1.2807813)(10.982814,1.0807812)
\pscircle[linewidth=0.04,dimen=outer](4.5,1.5807812){0.3}
\pscircle[linewidth=0.04,dimen=outer](5.5,1.5807812){0.3}
\pscircle[linewidth=0.04,dimen=outer](7.5,1.5807812){0.3}
\pscircle[linewidth=0.04,dimen=outer](10.5,1.5807812){0.3}
\psline[linewidth=0.04](4.0,-2.5192187)(4.0,-2.7192187)(8.0,-2.7192187)(8.0,-2.5192187)
\psline[linewidth=0.04cm](5.0,-2.5192187)(5.0,-2.7192187)
\psline[linewidth=0.04cm](6.0,-2.5192187)(6.0,-2.7192187)
\psline[linewidth=0.04cm](7.0,-2.5192187)(7.0,-2.7192187)
\pscircle[linewidth=0.04,dimen=outer,fillstyle=solid,fillcolor=black](4.5,-2.2192187){0.3}
\pscircle[linewidth=0.04,dimen=outer,fillstyle=solid,fillcolor=black](4.5,-1.4192188){0.3}
\pscircle[linewidth=0.04,dimen=outer,fillstyle=solid,fillcolor=black](6.5,-2.2192187){0.3}
\pscircle[linewidth=0.04,dimen=outer,fillstyle=solid,fillcolor=black](7.5,-2.2192187){0.3}
\pscircle[linewidth=0.04,dimen=outer,fillstyle=solid,fillcolor=black](7.5,-1.4192188){0.3}
\psarc[linewidth=0.02,arrowsize=0.05291667cm 2.0,arrowlength=1.4,arrowinset=0.4]{<-}(2.0,1.6807812){0.6}{30.963757}{149.03624}
\psarc[linewidth=0.02,arrowsize=0.05291667cm 2.0,arrowlength=1.4,arrowinset=0.4]{<-}(6.0,1.6807812){0.6}{30.963757}{149.03624}
\psarc[linewidth=0.02,arrowsize=0.05291667cm 2.0,arrowlength=1.4,arrowinset=0.4]{<-}(8.0,1.6807812){0.6}{30.963757}{149.03624}
\psarc[linewidth=0.02,arrowsize=0.05291667cm 2.0,arrowlength=1.4,arrowinset=0.4]{<-}(11.0,1.6807812){0.6}{30.963757}{149.03624}
\psarc[linewidth=0.02,arrowsize=0.05291667cm 2.0,arrowlength=1.4,arrowinset=0.4]{->}(10.0,1.6807812){0.6}{30.963757}{149.03624}
\usefont{T1}{ptm}{m}{n}
\rput(9.871407,2.5907812){$Q_{(5)}$}
\usefont{T1}{ptm}{m}{n}
\rput(11.051406,2.5907812){$P_{(5)}$}
\psarc[linewidth=0.02,arrowsize=0.05291667cm 2.0,arrowlength=1.4,arrowinset=0.4]{->}(1.0,1.6807812){0.6}{30.963757}{149.03624}
\psarc[linewidth=0.02,arrowsize=0.05291667cm 2.0,arrowlength=1.4,arrowinset=0.4]{->}(4.0,1.6807812){0.6}{30.963757}{149.03624}
\psarc[linewidth=0.02,arrowsize=0.05291667cm 2.0,arrowlength=1.4,arrowinset=0.4]{->}(7.0,1.6807812){0.6}{30.963757}{149.03624}
\usefont{T1}{ptm}{m}{n}
\rput(0.87140626,2.5907812){$Q_{(1)}$}
\usefont{T1}{ptm}{m}{n}
\rput(3.8714063,2.5907812){$Q_{(2)}$}
\usefont{T1}{ptm}{m}{n}
\rput(6.871406,2.5907812){$Q_{(4)}$}
\usefont{T1}{ptm}{m}{n}
\rput(2.0514061,2.5907812){$P_{(1)}$}
\usefont{T1}{ptm}{m}{n}
\rput(6.0514064,2.5907812){$P_{(3)}$}
\usefont{T1}{ptm}{m}{n}
\rput(8.051406,2.5907812){$P_{(4)}$}
\psarc[linewidth=0.02,arrowsize=0.05291667cm 2.0,arrowlength=1.4,arrowinset=0.4]{<-}(8.0,-0.51921874){0.6}{30.963757}{149.03624}
\psarc[linewidth=0.02,arrowsize=0.05291667cm 2.0,arrowlength=1.4,arrowinset=0.4]{<-}(7.0,-1.3192188){0.6}{30.963757}{149.03624}
\psarc[linewidth=0.02,arrowsize=0.05291667cm 2.0,arrowlength=1.4,arrowinset=0.4]{<-}(5.0,-0.51921874){0.6}{30.963757}{149.03624}
\psarc[linewidth=0.02,arrowsize=0.05291667cm 2.0,arrowlength=1.4,arrowinset=0.4]{<-}(4.0,-1.3192188){0.6}{30.963757}{149.03624}
\psarc[linewidth=0.02,arrowsize=0.05291667cm 2.0,arrowlength=1.4,arrowinset=0.4]{->}(6.0,-1.3192188){0.6}{30.963757}{149.03624}
\psarc[linewidth=0.02,arrowsize=0.05291667cm 2.0,arrowlength=1.4,arrowinset=0.4]{->}(7.0,-0.51921874){0.6}{30.963757}{149.03624}
\psarc[linewidth=0.02,arrowsize=0.05291667cm 2.0,arrowlength=1.4,arrowinset=0.4]{->}(4.0,-0.51921874){0.6}{30.963757}{149.03624}
\psarc[linewidth=0.02,arrowsize=0.05291667cm 2.0,arrowlength=1.4,arrowinset=0.4]{->}(8.0,-1.3192188){0.6}{30.963757}{149.03624}
\usefont{T1}{ptm}{m}{n}
\rput(4.0714064,-0.40921876){$Q_{(1)}$}
\usefont{T1}{ptm}{m}{n}
\rput(5.0714064,0.39078125){$Q_{(2)}$}
\usefont{T1}{ptm}{m}{n}
\rput(7.0714064,-0.40921876){$Q_{(4)}$}
\usefont{T1}{ptm}{m}{n}
\rput(8.071406,0.39078125){$Q_{(5)}$}
\usefont{T1}{ptm}{m}{n}
\rput(4.0514064,0.39078125){$P_{(1)}$}
\usefont{T1}{ptm}{m}{n}
\rput(6.0514064,-0.40921876){$P_{(3)}$}
\usefont{T1}{ptm}{m}{n}
\rput(7.0514064,0.39078125){$P_{(4)}$}
\usefont{T1}{ptm}{m}{n}
\rput(8.051406,-0.40921876){$P_{(5)}$}
\end{pspicture} 
}
\end{center}
\caption{Mapping of a simple exclusion process (above) to ZRP (below) for $n=5$.}
\label{fig:mapping}
\end{figure}
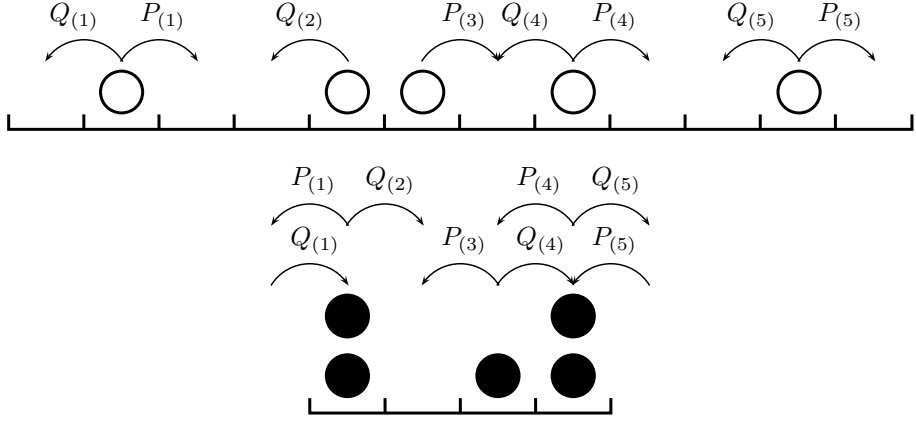

First notice that the dynamics of a multiple $n$-shock, which has been shown to be an $n$-particle simple exclusion dynamics with hopping rates $P_{(k)}$ (right) and $Q_{(k)}$ (left) for particle $k$, can be mapped to an $n-1$ site open zero range process (ZRP), see Figure \ref{fig:mapping}.%
The inter-particle distances in the exclusion process are interpreted as occupation numbers of the corresponding open ZRP. In this ZRP the left and right hop rates from site $k$ ($k=1,2,\cdots,n-1$) are independent of the occupation number and are equal to $Q_{(k+1)}$ and $P_{(k)}$ respectively. 
The rate of particle injection on the first and last ($n-1^\text{st}$) site of the system is $Q_{(1)}$ and $P_{(n)}$. 

Due to the attraction of ``micro-shocks'' in the ASEP, this ZRP has a well-defined stationary measure, which has a product structure with site-dependent fugacities \cite{Levine2005}. The mean velocity of the shock in the ASEP maps to the stationary current in the ZRP. In \cite{Levine2005} the stationary current is calculated for open zero range processes with constant left and right hopping rates in the bulk. Their formula (eq.~15 in that paper) for the current $J$ can easily be generalized for inhomogeneous cases as
\begin{equation}\label{vshock}
 J= \frac{\prod\limits_{k=1}^{n}P_{(k)} - \prod\limits_{k=1}^{n}Q_{(k)}}{\sum\limits_{k=1}^n \prod\limits_{\ell=1}^{k-1}P_{(\ell)} \prod\limits_{\ell=k+1}^{n}Q_{(\ell)}}.
\end{equation}
Using \eqref{jumprates} we write the denominator of \eqref{vshock} as
\begin{equation}\label{denominator}
 \sum_{k=1}^nq^{k-1}\prod_{\ell=1}^{k-1}\frac{\vr_{(\ell)}}{\vr_{(\ell-1)}}\cdot q^{n-k}\prod_{\ell=k+1}^n\frac{1-\vr_{(\ell-1)}}{1-\vr_{(\ell)}}
=\frac{q^{n-1}}{\vr_{(0)}(1-\vr_{(n)})}\sum_{k=1}^n\vr_{(k-1)}(1-\vr_{(k)}).
\end{equation}
Define now
\[
A=\sum_{k=1}^n\vr_{(k-1)}(1-\vr_{(k)}),\qquad B=\sum_{k=1}^n(1-\vr_{(k-1)})\vr_{(k)},
\]
for which we have
\[
p\cdot A=q\cdot B\qquad\text{and}\qquad B-A=\vr_{(n)}-\vr_{(0)}.
\]
This can be solved to
\[
A=\frac{\vr_{(n)}-\vr_{(0)}}{p-q}\cdot q
\]
and hence \eqref{denominator} equals
\[
\frac{\vr_{(n)}-\vr_{(0)}}{\vr_{(0)}(1-\vr_{(n)})}\cdot\frac{q^n}{p-q}.
\]
Inserting this into \eqref{vshock} one obtains
\begin{equation}\label{result}
J=\frac{\frac{\vr_{(n)}}{\vr_{(0)}}\cdot q^n-\frac{1-\vr_{(0)}}{1-\vr_{(n)}}\cdot q^n}{\frac{\vr_{(n)}-\vr_{(0)}}{\vr_{(0)}(1-\vr_{(n)})}\cdot\frac{q^n}{p-q}}=(p-q)\cdot\frac{\vr_{(n)}(1-\vr_{(n)})-\vr_{(0)}(1-\vr_{(0)})}{\vr_{(n)}-\vr_{(0)}}.
\end{equation}
One can see that in the rhs.\ of \eqref{result} the usual Rankine-Hugoniot velocity is recovered as expected.

\bibliographystyle{unsrt}
\bibliography{refsmarton}

\end{document}